\documentclass[11pt,letterpaper]{amsart}

\usepackage{enumitem}
\usepackage{amsfonts}
\usepackage{amssymb}
\usepackage{amsmath}
\usepackage{amsthm}
\usepackage{mathtools}
\usepackage{amstext}
\usepackage{tikz-cd}
\usepackage[left=1.3in,right=1.3in,top=1.2in, bottom=1.3in]{geometry}
\usepackage{graphicx}
\definecolor{cit}{HTML}{117733}
\definecolor{lin}{HTML}{0003d6}
\usepackage{hyperref}
\hypersetup{bookmarksnumbered,colorlinks=true,linkcolor=lin,citecolor=cit,urlcolor=black}

\usepackage[section]{placeins}
\usepackage{caption}

\DeclareMathOperator{\PGL}{PGL}
\DeclareMathOperator{\kum}{Kum}

\DeclareMathOperator{\Ind}{Ind}
\DeclareMathOperator{\Bir}{Bir}
\DeclareMathOperator{\Aut}{Aut}
\DeclareMathOperator{\PsAut}{PsAut}
\newcommand{\mzb}[1]{\overline{M}_{0,#1}}

\DeclareMathOperator{\proj}{Proj}

\DeclareMathOperator{\bl}{Bl}

\DeclareMathOperator{\Nef}{Nef}

\DeclareMathOperator{\Effb}{\overline{Eff}}

\DeclareMathOperator{\Mov}{Mov}
\DeclareMathOperator{\Pic}{Pic}
\DeclareMathOperator{\NS}{NS}

\DeclareMathOperator{\rank}{rank}
\DeclareMathOperator{\Span}{span}

\DeclareMathOperator{\id}{id}

\DeclareMathOperator{\spec}{Spec}

\newcommand{\CC}{\mathbb{C}}
\newcommand{\Z}{\mathbb{Z}}
\newcommand{\F}{\mathbb{F}}
\newcommand{\R}{\mathbb{R}}
\newcommand{\PP}{\mathbb{P}}

\newcommand{\Q}{\mathbb{Q}}
\newcommand{\oO}{\mathcal{O}}

\newcommand{\e}{\equiv}
\newcommand{\m}{\pmod} 

\newcommand{\vv}[1]{\lvert#1\rvert}
\newcommand{\vr}[1]{\langle#1\rangle}
\newcommand{\ra}{\to}

\newcommand{\mt}{\mapsto}

\newcommand{\ik}{^{-1}}

\newcommand{\txt}[1]{\text{ #1}} 
\theoremstyle{plain}
\newtheorem{theorem}{Theorem}[section]
\newtheorem{defiThm}[theorem]{Definition-Theorem}
\newtheorem{corollary}[theorem]{Corollary}
\newtheorem{prop}[theorem]{Proposition}
\newtheorem{lemma}[theorem]{Lemma}
\newtheorem{question}[theorem]{Question}

\theoremstyle{definition}
\newtheorem{defi}[theorem]{Definition}

\newtheorem{remark}[theorem]{Remark}

\newcommand{\pf}{{\em Proof}}
\newcommand{\pfof}[1]{{\em Proof of #1}}

\begin{document}
\title[Birational geometry of blow-ups of $\PP^n$ along points and lines]{Birational geometry of blow-ups of projective spaces along points and lines}
\author{Zhuang He}\email{he.zhu@northeastern.edu}\address{Department of Mathematics, Northeastern University, 360 Huntington Ave, Boston MA}
\author{Lei Yang}\email{yang.lei1@northeastern.edu}\address{Department of Mathematics, Northeastern University, 360 Huntington Ave, Boston MA}
\date{}

\begin{abstract}
	Consider the blow-up $X$ of $\PP^3$ at $6$ points in very general position and the $15$ lines through the $6$ points. We construct an infinite-order pseudo-automorphism $\phi_X$ on $X$, induced by the complete linear system of a divisor of degree $13$. The effective cone of $X$ has infinitely many extremal rays and hence, $X$ is not a Mori Dream Space. The threefold $X$ has a unique anticanonical section which is a Jacobian K3 Kummer surface $S$ of Picard number 17. The restriction of $\phi_X$ on $S$ realizes one of Keum's 192 infinite-order automorphisms of Jacobian K3 Kummer surfaces. In general, we show the blow-up of $\PP^n$ ($n\geq 3$) at $(n+3)$ very general points and certain $9$ lines through them is not Mori Dream, with infinitely many extremal effective divisors. As an application, for $n\geq 7$, the blow-up of $\mzb{n}$ at a very general point has infinitely many extremal effective divisors. 
\end{abstract}

\maketitle

\section{Introduction}
\label{Intro}
\begin{figure}[h]
	\centering
	\includegraphics[height=3.8cm,clip]{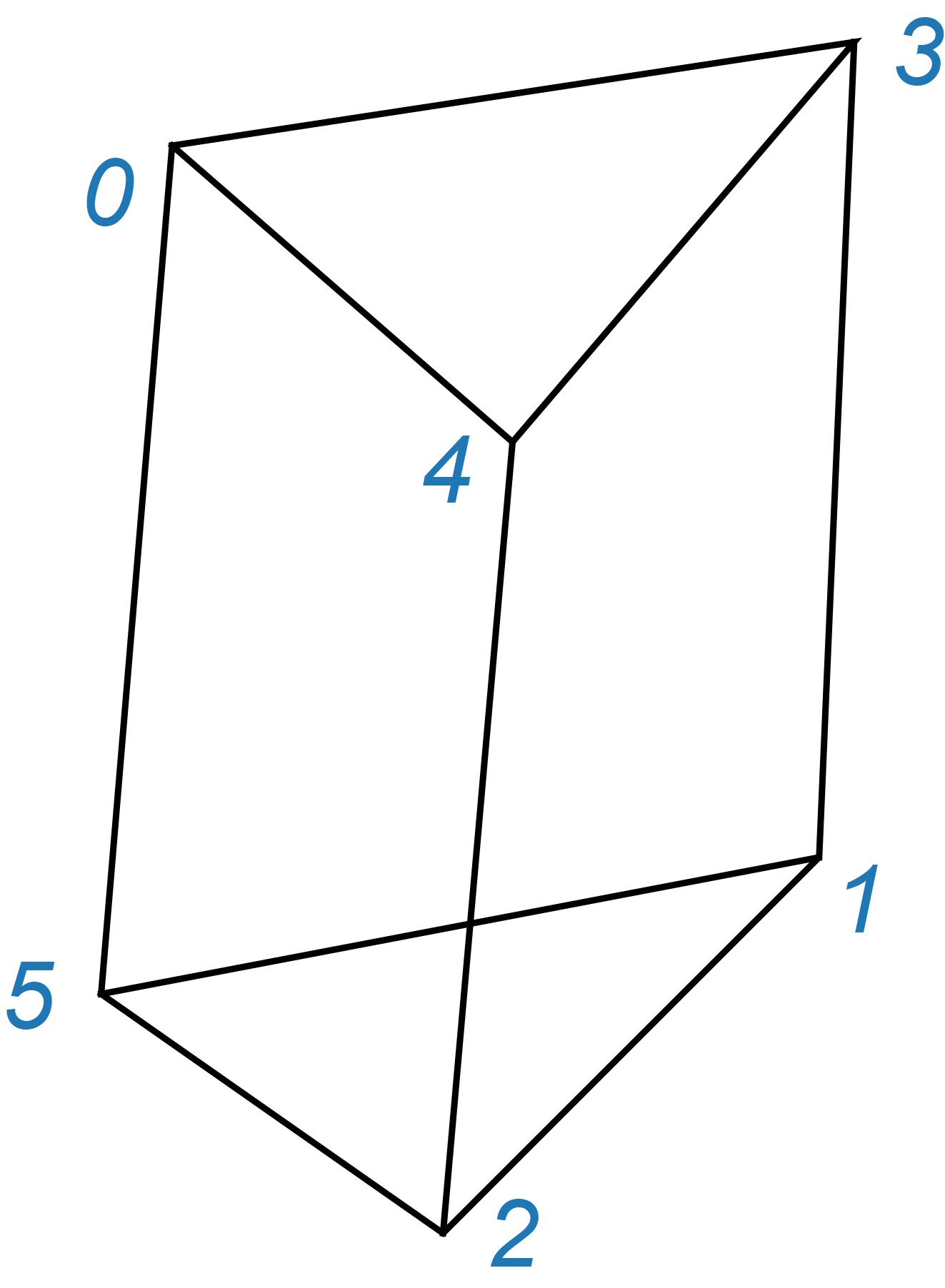}
	\caption{The $9$ lines blown-up in $\PP^3$}
	\label{fig:prism}
\end{figure}

{\let\thefootnote\relax\footnote{{This project was partially supported by the NSF grant DMS-1701752.}}}
We consider the blow-ups of the projective space $\PP^3$ at points and lines. We work over the complex numbers.
Define:
\begin{itemize}
	\item $u:Y\ra \PP^3$ to be the successive blow-up of $\PP^3$ at $6$ points $p_0,\cdots, p_5$ in (very) general position, and the proper transforms of the $9$ lines $\overline{p_i p_j}$ labeled by 
\[(ij)\in\mathcal{I}=\{03,04,34,12,15,25,05,13,24\};\]
\item $v:X\ra \PP^3$ to be the successive blow-up of $\PP^3$ at $p_0,\cdots,p_5$ and the proper transforms of all the $15$ lines $\overline{p_i p_j}$.
\end{itemize}
The configuration of the $9$ lines blown-up to get $Y$ is best shown in Figure \ref{fig:prism}.

A birational map $f:Y\dashrightarrow Y$ is called a pseudo-automorphism if there are open sets $U$ and $V$ in $Y$ such that $f:U\ra V$ is an isomorphism, and the codimensions of $Y\backslash U$ and $Y\backslash V$ in $Y$ are at least $2$.
The main construction of this paper is an infinite-order pseudo-automorphism $\phi$ of $Y$. 
Let $E_i$ and  $E_{ij}$ be the exceptional divisors of the blow-ups $Y$ and $X$ over the points $p_i$ and lines $\overline{p_i p_j}$. Let $H:=u^*\oO_{\PP^3}(1)$. Then the Picard group of $Y$ is freely generated by $H, E_i$ and $E_{ij}$.
Consider the following divisor class over $Y$ (and $X$):
\begin{align}\label{Dclass}
	\begin{split}
	D:=\quad &  13H-7(E_1+E_2+E_5)-5(E_0+E_3+E_4)\\
	&-3(E_{03}+E_{04}+E_{34})
	-4(E_{05}+E_{13}+E_{24})-(E_{12}+E_{15}+E_{25}).
	\end{split}
\end{align}
We point out that the divisor class $D$ is fixed by the action of $\mathcal{S}_3$ on the ordered pairs of the six points $\{ (5,0),(1,3),(2,4)\}$. That is, the action which permutes the vertical edges of the prism in Figure \ref{fig:prism} while keeping their directions.

Recall that Mori Dream Spaces are introduced by \cite{Hu2000} (see Section \ref{prebm} for definition and properties). The effective cone of a Mori Dream Space is rational polyhedral, with a chamber decomposition which determines its birational geometry.
We have the following results:

\begin{theorem}\label{main}
	For very general six points $p_0,\cdots, p_5$:
	\begin{enumerate}
		\item The linear system $\vv{D}$ has dimension $3$; hence it determines a map $\phi_D:Y\dashrightarrow \PP^3$. 
		\item There exist $6$ points $q_0,\cdots, q_5$ in the target copy of $\PP^3$ which are projectively equivalent to $p_0,\cdots, p_5$. That is, there exists $M\in \PGL(4)$ such that $M p_i=q_i$ for $i=0,\cdots, 5$. Blowing up the $6$ points $q_i$ and the corresponding $9$ lines $\overline{q_i q_j}$ for $(ij)\in \mathcal{I}$ induces a pseudo-automorphism $\phi:Y\dashrightarrow Y$. Blowing up $q_i$ and all the $15$ lines $\overline{q_i q_j}$ induces a pseudo-automorphism $\phi_X:X\dashrightarrow X$.
		\item The pseudo-automorphisms $\phi$ and $\phi_X$ are of infinite order.
		\item The effective cone $\Effb(Y)$ of $Y$ has infinitely many extremal rays (see Theorem \ref{inforder}). Hence $\Effb(Y)$ is not rational polyhedral, and $Y$ is not a Mori Dream Space. The same results hold for $X$.
	\end{enumerate}	
\end{theorem}

We note that the question whether $X$ is a Mori Dream Space was proposed by John Ottem.

The divisor class $D$ and the pseudo-automorphism $\phi$ are related to Keum's automorphisms of Jacobian K3 Kummer surfaces. A Kummer surface $\kum(A)$ is the quotient of an abelian surface $A$ under the involution $\iota:A\ra A, a\mt -a$. The set of order-$2$ points on $A$, denoted by $A[2]$, has $16$ elements. The surface $\kum(A)$ is singular with $16$ nodes over $A[2]$. The minimal desingularization of $\kum(A)$ is a K3 surface $K(A)$, which we refer to as the K3 Kummer surface associated with $A$. We say $K(A)$ is of Jacobian type if $A\cong J(C)$ is the Jacobian variety of a smooth genus-$2$ curve $C$.

In our context, the key fact is that $X$ has a unique anticanonical section $S$ which is a smooth K3 Kummer surface of Jacobian type, with Picard rank $\rho(S)=17$, for very general six points $p_i$ in $\PP^3$. Keum \cite{Keum1997} first constructed $192$ infinite-order automorphisms of a Jacobian K3 Kummer surface $S$ of Picard rank $17$, each associated with one of $192$ Weber Hexads, which are certain $6$-element subsets of $A[2]$
.  If we denote by $\PsAut(X)$ the group of pseudo-automorphisms of $X$, then restricting to $S$ induces a group homomorphism $s:\PsAut(X) \ra \Aut(S)$, for the reason that $S$ is K3 and is the unique anticanonical section of $X$. In fact, $s(\phi_X)$ is one of these $192$ automorphisms:
\begin{theorem}
	\label{k3KeumAuto}
For very general six points $p_0,\cdots, p_5$:
	\begin{enumerate}
		\item $X$ has a unique anticanonical section $S$, which is a Jacobian K3 Kummer surface with $\rho(S)=17$.
		\item The restriction of $\phi_X$ to $S$ equals Keum's automorphism $\kappa:S\ra S$ associated with the Weber Hexad $\mathcal{H}=\{1,2,5,12,14,23\}$ (see Section \ref{k192}).
		\item The inverse $\phi_X\ik$ (and $\phi\ik$) is induced by the complete linear system of $D'$ where 
\begin{align*}\label{Dprimeclass}
	\begin{split}
	D':=\quad &  13H-5(E_1+E_2+E_5)-7(E_0+E_3+E_4)\\
	&-(E_{03}+E_{04}+E_{34})-4(E_{05}+E_{13}+E_{24})-3(E_{12}+E_{15}+E_{25}).
	\end{split}
\end{align*}
	\end{enumerate}	
\end{theorem}

In particular, there are $60$ such configurations $\mathcal{I}$, each deciding a pair of pseudo-automorphisms inverse to each other. In total we have $120$ such pseudo-automorphisms. Their restrictions to $S$ are exactly $120$ out of the $192$  Keum's automorphisms.

We consider the birational automorphism $\psi:\PP^3\dashrightarrow\PP^3$ induced by $\vv{D}$. It turns out that $\psi$ contracts exactly $9$ distinct irreducible rational quartics $Q_\alpha$, indexed by $\alpha\in \mathcal{A}:=\{0,3,4,12,15,25,05,13,24\}$. We refer to Section \ref{9Q} for their divisor classes. Here we summarize the key features of these quartics and the divisor class $D$:

\begin{enumerate}
	\item Each $Q_\alpha$ is unique in its divisor class when considered over $Y$ or $X$ (see Theorem \ref{pencilsQ}).
	\item $\phi$ maps $Q_\alpha$ birationally onto the exceptional divisor $E_\alpha$ (Section \ref{imageQ} and Theorem \ref{psdautoY}).
	\item The restrictions of $Q_\alpha$ to $S$ are all $(-2)$-curves (see Theorem \ref{pencilsQ}).
	\item The restriction of $D$ to $S$ equals $\kappa ^* H_S$, where $H_S$ is a hyperplane class of $S$ (see Formula (\ref{rel4}) and Theorem \ref{pencilsQ}).
\end{enumerate}

The birational map $\psi$ also has a surprising interaction with the structure theory of $\Bir(\PP^3)$, the birational automorphism group of $\PP^3$. A classical result by Max Noether and Castelnuovo \cite{Castelnuovo1901} says that $\Bir(\PP^2)$ is generated by $\PGL(3)$ and the standard Cremona transformation $\sigma_2:[x:y:z]\mt [1/x:1/y:1/z]$. The analogue is false for $n\geq 3$, where $\Bir(\PP^n)$ is strictly larger than the subgroup  $G_n:=\vr{\PGL(n+1),\sigma_{n}}$ \cite{Hudson1927}\cite{Pan1999}, and $\sigma_n$ is the standard Cremona transformation of $\PP^n$. One of the interesting subsets of the large group $\Bir(\PP^n)$ is the set $H_n$ of all $f\in \Bir(\PP^n)$ which only contracts rational hypersurfaces. It is known that $G_n\subset H_n$ (see \cite[\S 1]{Blanc2014}). On the other direction, \cite{Blanc2014} proved that $G_n\neq H_n$ when $n\geq 3$ is odd over any field $\mathbf{k}$, by giving examples of monomial birational maps which only contract rational hypersurfaces but are not in $G_n$ when $n$ is odd. They further gave a criterion \cite[Thm. 1.4]{Blanc2014} (see Theorem \ref{BH}) characterizing elements in $G_n$, which we can apply to  $\psi$ and find: 
\begin{theorem}
	Over $\CC$, $\psi\in H_3$ but $\psi\not\in G_3$.
	\label{psinotG3}
\end{theorem}

In general, we consider the successive blow-up of $\PP^n$ at points and lines.   By \cite{Mukai2001} and \cite{Castravet2006}, the blow-up of $\PP^n$ at $r$ very general points $p_1,\cdots, p_r$ is a Mori Dream Space if and only if its effective cone is rational polyhedral, which in turn is equivalent to 
	\[\frac{1}{n+1}+\frac{1}{r-n-1}>\frac{1}{2}.\]

In particular, the last inequality translates to $r\leq 8$ for $n=2,4$, $r\leq 7$ for $n=3$ and $r\leq n+3$ for $n\geq 5$. There are many further results on the birational geometry of $\PP^n$ blown-up at points and lines, including \cite{Sturmfels2010,Araujo2015,Brambilla2016,Coskun2016,Araujo2017,Dumitrescu2017,Casagrande2019,Pintye2019}. 

\begin{question}\label{n3Q}
	Let $X'$ be the blow-up of $\PP^n$ ($n\geq 3$) at $(n+3)$ points in general position and certain lines through the $(n+3)$ points. For what configuration of the lines is $\Effb(X')$ rational polyhedral?
\end{question}
We shed light upon Question \ref{n3Q} by showing that as soon as we blow up $9$ lines in a very special configuration, the effective cone of the blow-up is no longer rational polyhedral. 

For $n\geq 3$, we define $Y_n$ to be the blow-up of $\PP^n$ at $(n+3)$ points in very general position and $9$ lines through six of them, such that when the six points are indexed by $0$ to $5$, the $9$ lines are labeled by $\mathcal{I}$ (see Figure \ref{fig:prism}). In particular, $Y_3=Y$ as defined above.
\begin{theorem}\label{main2}
	For each $n\geq 4$ there is a small $\Q$-factorial modification (SQM) (see Section \ref{prebm} for definition) $\tilde{Y}_n$ of $Y_n$ such that $\tilde{Y}_n$ is a $\PP^1$-bundle over $Y_{n-1}$.
	For $n\geq 3$, $\Effb(Y_n)$ has infinitely many extremal rays. Hence $Y_n$ are not Mori Dream for $n\geq 3$.	
\end{theorem}

Let $\overline{M}_{g,n}$ be the Deligne-Mumford compactification of the moduli space of stable curves of genus $g$ with $n$ marked points. One of the questions of the birational geometry of $\overline{M}_{0,n}$ is to determine whether they are Mori Dream Spaces. Castravet and Tevelev \cite{Castravet2015} first proved that $\mzb{n}$ is not a Mori Dream Space for $n>133$, which was later improved by \cite{Gonzalez2016} and \cite{Hausen2018} to $n\geq 10$. On the other hand, for $n\leq 6$, $\overline{M}_{0,n}$ are of Fano type, and hence Mori Dream Spaces.
\begin{question}
	Is $\mzb{n}$ Mori Dream for $n=7,8$ and $9$?
\end{question}
Recall Kapranov's blow-up construction of $\mzb{n}$ \cite{Kapranov1993} which realizes $\mzb{n}$ as the successive blow-up of linear subspaces of $\PP^{n-3}$ of codimensions at least $2$ passing through points among $(n-1)$ points in linearly general position, in increasing order. Now the blow-up at one more point of the Kapranov's blow-up model of $\mzb{n}$ factors through $Y_{n-3}$ if $n\geq 7$. We have the following result.
\begin{theorem}
	\label{mzbblowup}
For $n\geq 7$, the effective cone of the blow-up of $\mzb{n}$ at a very general point has infinitely many extremal rays. Hence the blow-up of $\mzb{n}$ at a very general point is not a Mori Dream Space.
\end{theorem}
We note that when $n=6$, the one-point blow-up of $\mzb{6}$ is not a blow-up of $Y_3=Y$, so Theorem \ref{mzbblowup} does not extend to $n=6$. The not Mori Dream Space part of Theorem \ref{mzbblowup} is new for $n=7,8$ and $9$, and for $n\geq 10$ it follows from that $\mzb{n}$ is not a Mori Dream Space and Okawa's result \cite{Okawa2016} (see Section \ref{preMDS}).

Structure of this paper: in Section \ref{prebm} we review generalities on birational maps, pseudo-isomorphisms and Mori Dream Spaces. Sections \ref{JKK} and \ref{uniquekk} show that $X$ has a unique anticanonical section $S$, which is a Jacobian K3 Kummer surface of Picard rank $17$ when the six points are very general. This proves Theorem \ref{k3KeumAuto} (1) (2). Section \ref{uniquekk} also identifies the Picard lattice of $S$ with the description which comes from the Jacobian Kummer structure.

Section \ref{k192} discusses the relations among the divisors over $S$, and then reviews Keum's $192$ automorphisms. 
Sections \ref{9Q} and \ref{dD} define the $9$ quartics $Q_\alpha$ and use them to build various sections of $D$. Section \ref{BIR} proves that the six points $\{q_i\}$ are projectively equivalent to $\{p_i\}$, and the rational map $\psi$ induced by $\vv{D}$ is birational. Section \ref{imageQ} shows that $\phi$ contracts none of the quartics $Q_\alpha$, and Section \ref{imageE} shows that $\phi$ does not contract any exceptional divisors. Thus we finish the proofs of Theorems \ref{main} and \ref{k3KeumAuto} in Theorems \ref{psdautoY}, \ref{inforder}, \ref{restrict} and Corollary \ref{psdautoX}. Section \ref{Cremona} relates $\psi$ to the birational automorphism group of $\PP^3$ and proves Theorem \ref{psinotG3}. Section \ref{SQMpn} is the application to the blow-up of $\PP^n$ and $\mzb{n}$, where we prove Theorem \ref{main2} and \ref{mzbblowup}.

{\bf Acknowledgement} We are grateful to our advisor Ana-Maria Castravet for suggesting this question, and for her constant support and countless discussions on this project. We thank J\'em\'ery Blanc, Igor Dolgachev, Frederic Han, Kalle Karu, Antonio Laface, Yucheng Liu, Emanuele Macr\`{i}, Karl Schwede and Jenia Tevelev for their helpful discussions and suggestions. We thank the referees for valuable comments. We used {\em Macaulay2} \cite{M2} and the {\em RationalMaps} package \cite{RationalMapsSource} to verify some of our results. Zhuang thanks University of Versailles for the hospitality during the visit when part of the paper was written. 

\section{Preliminaries on Birational Maps and Mori Dream Spaces}\label{prebm}
\subsection{Birational maps and pseudo-automorphisms}\label{prebpa} We work over $\CC$. Let $X$, $Y$ be normal projective varieties. We say $f$ is a pseudo-isomorphism (See \cite[2.2]{Cantat2019}) if $f$ is birational and there exist Zariski open subsets $U\subset X$ and $V\subset Y$ such that (1) $f_{\mid_U} :U\ra V$ is an isomorphism and (2) $X\backslash U$ and $Y\backslash V$ have codimension at least $2$. For the case $X=Y$ we say $f$ is a pseudo-automorphism of $X$.

The indeterminacy locus $\Ind(f)$ of $f$ is defined to be $X-U_{0}$ where $U_{0}$ is the largest open subset of $X$ on which $f$ is defined. When $X$ and $Y$ are normal and projective, $\Ind(f)$ of $f$ has codimension $\geq 2$. Then we can define the image $f(Z)$ of a codimension $1$ subvariety $Z\subset X$ as the Zariski closure of $f(Z\backslash \Ind(f))$. We say $f$ contracts $Z$ if the codimension of $f(Z)$ in $Y$ is at least $2$.
We recall the following basic fact (see \cite[Prop. 2.4]{Cantat2019} for a proof):
\begin{lemma}\label{psdhyp}
	Let $X$, $Y$ be normal projective varieties and let $f:X\dashrightarrow Y$ be a birational map. Then $f$ is a pseudo-isomorphism if and only if neither $f$ nor $f\ik$ contracts any divisors.
\end{lemma}

Given a birational map $f:X\dashrightarrow Y$, the Jacobian determinant $\det J(f)(x)$ of $f$ at a point $x\in X$ can be defined as the determinant of $d f_x$ in some local coordinates. The value $\det J(f)(x)$ depends on the local coordinates, but whether $\det J(f)(x)=0$ does not. Furthermore, $\det J(f)(x)\neq 0$ if and only if $f$ is locally an isomorphism at $x$, or equivalently, $f$ is \'{e}tale at $x$.  Therefore we can define the exceptional set of $f$ as the subset of $X$ where $f$ is not defined or locally not an isomorphism.

In the special case when $f:\PP^n \dashrightarrow \PP^n$ is a birational automorphism of $\PP^n$, $f$ is defined by $[f_0:\cdots:f_n]$ for homogeneous degree-$d$ polynomials $f_i\in \CC[x_0,\cdots,x_n]$, with $\gcd(f_0,\cdots, f_n)=1$. In this case we have $\det J(f)=\det \left( \partial f_i/\partial x_j\right)_{0\leq i,j\leq n}$.
Since $f$ is birational, we must have $\det J(f)\not\equiv 0$ is a nonzero polynomial of degree at most $m=(d-1)(n+1)$. 
When $m\geq 1$, $\det J(f)$ defines the exceptional set of $f$ \cite[7.1.4]{Dolgachev2012}, which is a hypersurface of degree at most $m$, and is the union of all the irreducible hypersurfaces contracted by $f$.

\subsection{Mori Dream Spaces} \label{preMDS}
A variety $X$ is $\Q$-factorial if for any Weil divisor $D$ on $X$, there exists some integer $m$ such that $mD$ is Cartier. For instance, smooth varieties are $\Q$-factorial.
A small $\Q$-factorial modification (SQM) of $X$ is a rational map $g:X\dashrightarrow X'$ such that $X'$ is $\Q$-factorial and $g$ is an isomorphism in codimension $1$.

By \cite{Hu2000}, A normal projective $\Q$-factorial variety $X$ is a Mori Dream Space (MDS) if
\begin{enumerate}
	\item $\Pic(X)$ is finitely generated;
	\item $\Nef(X)$ is spanned by finitely many semiample divisors.
	\item There are finitely many SQMs $g_i:X\dashrightarrow X_i$ such that each $X_i$ satisfies (1) and (2) above, and the movable cone $\Mov(X)$ is the union of $g_i^* \Nef(X_i)$. 
\end{enumerate}

By definition, if $X$ is a Mori Dream Space, then any SQM $X_i$ of $X$ is a Mori Dream Space. Later we will use the following result by Okawa \cite{Okawa2016}. Suppose $X$ and $Y$ are normal, projective, $\Q$-factorial varieties and $f:X\ra Y$ is a surjective morphism. If $X$ is a Mori Dream Space, then $Y$ is also a Mori Dream Space.

\section{Preliminary on K3 Kummer Surfaces}\label{JKK}
Kummer surfaces are classically defined as singular quartics in $\PP^3$ with $16$ nodes. Here we adopt the following definition:
\begin{defi}\label{singK}
	Let $A$ be an algebraic abelian surface. Then the (singular) {\em Kummer surface} $\kum(A)$ associated with $A$ is the quotient $A/\iota$ where $\iota:A\ra A, a\mt -a$ is the involution.
\end{defi}
An abelian surface $A$ has exactly $16$ order-$2$ points, which form a group $A[2]\cong \Z_2^{\oplus 4}$. Therefore $\kum(A)$ is a singular surface with exactly $16$ nodes. Now we have $f:A\ra \kum(A)$ is the double cover branched at $A[2]$. Let $\pi:S\ra \kum(A)$ be the minimal de-singularization of $\kum(A)$. Then $S$ is a smooth K3 surface. We call $S$ the {\em K3 Kummer surface} associated with $A$.

Later we will always identify the $16$ nodes in $\kum(A)$ with $A[2]$. Blowing-up $A$ at $A[2]$ gives us a smooth surface $\tilde{S}$.  Denote by $\pi'$ the blow-up $\pi':\tilde{S}\ra A$. Then there is a double cover $f':\tilde{S}\ra \tilde{S}/\tilde{\iota}\cong S$, where $\tilde{\iota}$ is an automorphism of $\tilde{S}$ lifting $\iota$,  with fixed locus being the $16$ exceptional divisors in $S$ over the $16$ nodes in $\kum(A)$. We have the following commutative diagram:
\begin{equation}
\begin{tikzcd}[]
	\tilde{S} \ar[d,"\pi'"] \ar[r,"f'"] & S \ar[d,"\pi"] \\   
	A \ar[r,"f"] & \kum(A)
\end{tikzcd}.
\label{cd}
\end{equation}

We recall the following key result on K3 Kummer surfaces (see \cite{Nikulin1975} and \cite[\S 14, Rem. 3.19]{Huybrechts2016}):
\begin{lemma}
	\label{NK16} 
	Let $S$ be a complex projective K3 surface. Then $S$ is isomorphic to a K3 Kummer surface if and only if there exist $16$ disjoint smooth rational curves on $S$.

Moreover, if there are $16$ disjoint smooth rational curves $C_i$ on the K3 surface $S$, then there exists an abelian surface $A$ and $\pi: S\ra \kum(A)$ such that $\pi$ is the de-singularization at the $16$ nodes with exceptional divisors $C_i$.
\end{lemma}

Recall that a curve $C$ on a K3 surface $S$ is called a $(-2)$-curve if $C$ is irreducible and $C^2=-2$. A $(-2)$-curve $C$ is necessarily isomorphic to $\PP^1$, and $h^0(S,\oO_S(C))=1$. We briefly review the lattice theory on K3 surfaces. If $S$ is a K3 surface, then $H^2(S,\Z)\cong E_8(-1)^{\oplus 2} \bigoplus U^{\oplus 3}$ is the K3 lattice  (see \cite[\S 1, Prop 3.5]{Huybrechts2016}). The Picard lattice $\Pic(S)$ and the N\'{e}ron-Severi lattice $\NS(S)$ coincide, and $T(S)=\NS(S)^\perp$ is the transcendental lattice of $S$. Consequently, numerically equivalence on $S$ is the same as linear equivalence, and any $(-2)$-curve $C$ is unique in its numerical class. 

We now return to $A[2]\cong \Z_2^{\oplus 4}$. Consider $A[2]$ as the affine 4-space over $\mathbb{F}_2$. Then there are exactly $30$ hyperplanes $\Gamma$ in $A[2]$, each containing $8$ elements. We have: 
\begin{lemma}\cite{Nikulin1975} (also see \cite[\S 14, Definition 3.13]{Huybrechts2016})
	Let $S$ be a K3 Kummer surface with exceptional divisors $N_\alpha$, $\alpha\in A[2]$. If $M\subset A[2]$ satisfies that $(1/2) \sum_{\alpha\in M} N_\alpha \in \NS(S)$, then $M$ is $\emptyset$, $A[2]$ or a hyperplane in $A[2]$.
	\label{NK2}
\end{lemma}

We say a Kummer surface $S$ associated with an abelian surface $A$ is of {\em Jacobian type} if $A=J(C)$ for some smooth genus $2$ curve $C$.
Next, we review that K3 Kummer surfaces are not Mori Dream Spaces.
Let $S$ be the K3 Kummer surface associated with an abelian surface $A$. The Picard ranks of $S$ and $A$ satisfy the relation \cite[\S 3, Rem 2.8]{Huybrechts2016} that 
\begin{align}
\rho(S)=\rho(A)+16.
	\label{add16}
\end{align}
Thus $17\leq \rho(S)\leq 20$.
By \cite{Artebani2010}, a K3 surface $S$ is a Mori Dream Space if and only if the effective cone $\Effb(S)$ of $S$ is rational polyhedral. By \cite{Pjateckiui-Sapiro1971}, a K3 surface with $\rho(S)\geq 3$ has rational polyhedral effective cone if and only if $\vv{\Aut(S)}<\infty$  (see \cite{Kovacs1994}). Therefore a K3 Kummer surface $S$ is a Mori Dream Space if and only if $\Aut(S)$ is finite. Now there is a complete classification of the Picard lattice of K3 surfaces with finite automorphism group for $\rho(S)\geq 17$  by \cite{Nikulin1983}, as well as a classification of the transcendental lattices $T(S)$ of K3 Kummer surfaces of $\rho(S)\geq 17$ \cite{Morrison1984}.  
A simple comparison shows that there is no compatibility in the two lists, noticing that $\NS(S)$ and $T(S)$ have the same determinant. Hence K3 Kummer surfaces are not Mori Dream Spaces. See \cite[Ex. 5.13]{Arzhantsev2014}.

\section{The Unique Anticanonical Section is K3 Kummer}\label{uniquekk}
In the rest of this paper, $X$ is the blow-up of $\PP^3$ at six points and the $15$ lines through them, with certain general position conditions on the points, which we will later specify.

\subsection{The K3 surface and $16$ disjoint lines}
 Let $H$ be the hyperplane class of $X$, and let $E_i$ and $E_{ij}$ be the exceptional divisors in $X$ over the points $p_i$ and lines $\overline{p_i p_j}$. Then $\Pic(X)$ is freely generated by $H, E_i$ and $E_{ij}$ over $\Z$.

\begin{lemma}
	Let $K_X$ be the canonical divisor of $X$. Then for the six points $p_0,\cdots, p_5$ in general position, $h^0(X, \oO_X(-K_X))=1$. The unique anticanonical section $S$ is a smooth K3 surface.
	\label{uniqueminusK}
\end{lemma}

We include a remark relating $S$ with the Weddle surfaces in spaces, which was kindly suggested by Professor Igor Dolgachev.
\begin{remark}
	Fix $p_0,\cdots, p_5$ in $\PP^3$ in linearly general position. The {\em Weddle surface} $W$ in $\PP^3$, first studied by Weddle \cite{Weddle1850}, is the locus of singular points of the net of singular quadrics passing through all those $p_i$. We note $W$ has nodal singularity at the six points $p_i$; $W$ contains the $15$ lines $\overline{p_i p_j}$ and the unique rational normal curve $R_0$ through the six points $p_i$.
	Here the image of the anticanonical section $S$ in $\PP^3$ is the Weddle surface $W$. 
	It is a classical result that $W$ is birational to the Jacobian Kummer surface $\kum(J(C_0))$ where $C_0$ is the genus $2$ curve double covering $R_0$ with branch locus $\{p_i\}$. Then our Propositions \ref{166} and \ref{S17} and Theorem \ref{kumaddition} are covered by \cite[\S 97, 98]{Hudson1905}. Also see \cite[\S 1]{Varley1986} for a modern treatment. In this section we provide different proofs using Nikulin's results on the lattice theory of Kummer surfaces, and properties of principally polarized abelian surfaces.
\end{remark}
\pfof{Lemma \ref{uniqueminusK}}. Here $-K_X\sim 4H-2\sum_i E_i -\sum_{ij} E_{ij}$. We place the points at $p_0=[1:A:B:C]$, with
$p_1=[1:0:0:0]$, $p_2=[0:1:0:0]$, $p_3=[0:0:1:0]$, $p_4=[0:0:0:1]$, $p_5=[1:1:1:1]$. We can assume $p_i$ are in linearly general position, and $\{1,A,B,C\}$ are distinct and nonzero.
Let $Z$ be the iterated blow-up of $\mathbb{P}^3$ along $p_1,\cdots,p_5$ and along the proper transforms of the $10$ lines $ \overline{p_ip_j}$. Then $Z$ is a Kapranov model of $\overline{M}_{0,6}$ (see Section \ref{SQMpn}). Direct calculation shows $h^0(Z,\oO_X(-K_Z))=5$.
 
 Let $\mathbb{P}^3=\proj \mathbb{C}[x,y,z,w]$. We construct five linearly independent sections of $\vv{-K_Z}$:
\[\begin{array}{ccc}
	f_1=(x-y)(z-w)xy, &\quad f_2=(x-y)(z-w)zw, & f_3=(x-z)(y-w)xz,\\
	f_4=(x-z)(y-w)yw, &\quad f_5=(x-w)(y-z)xw.
\end{array}\]
Then each $f_i$ is a quartic in $\PP^3$, vanishes at least twice on $p_1,\cdots,p_5$ and vanishes on $\overline{p_ip_j}$ for $1\leq i<j\leq 5$. It is easy to verify that  $f_i$ are linearly independent. Hence the proper transforms in $Z$ of $f_1,\cdots,f_5$ span $H^0(Z,\oO_Z(-K_Z))$.

First, we show the uniqueness of the anticanonical section of $X$. Suppose $S$ is an anticanonical section of $X$, with $S'$ its image in $\PP^3$. Let $f$ be a polynomial defining $S'$.  Then the proper transform of $S'$ in $Z$ is a section of $-K_Z$. Hence $f=af_1+bf_2+cf_3+df_4+ef_5$ for some constants $a,b,c,d$ and $e$ in $\mathbb{C}$. 
That $f$ vanishes on $\overline{p_0p_1}$ and $\overline{p_0p_2}$ implies $f([1,At,Bt,Ct])=0$ and $f([t,(A-1)t+1,Bt,Ct])=0$. Equivalently, we have the following:
\begin{align*}
	&ABC((C-B)b + (C-A)d)=0,\\
	& A^2( C-B)a -B C (C-B) b + B^2 (C-A) c - A C (C-A ) d + C^2(B-A) e=0,\\
	& A (C-B)a + B( C-A )c + C(B-A )e=0,\\
& ( A-2 ) ( A-1) (C-B) a+ (A-2) B (C-B) C b \\
& + (B-1) B ( 1-A + C) c + ( A-1) ( B-1) C (1-A+C) d+ ( 1-A+B) ( C-1) Ce=0,\\
&(2 A-3) (C-B)a+B C(C-B) b- B(B-1) c + C(B-1) (2-2 A + C)d - C(C - 1) e=0,\\
&(C-B) a -C(B- 1)d=0,
\end{align*}
which gives a matrix $M_{p_0}$ such that $M_{p_0}v=0$ for the vector $v=[a, b, c, d, e]^T$.
When $p_0=[1:2:3:4]$, we can directly compute that 
$\rank(M_{p_0})=4$. This implies that for $p_0$ general, $\rank(M_{p_0})\geq 4$.  As a result, for $p_0$ general, $S'$ must be unique if it exists. Therefore if $S$ exists, it must be unique.

For existence, let $g_1=(B - 1)C$, $g_2=(A-C)$, $g_3=-(A-1)C$, $g_4=-(B-C)$ and $g_5=(AB-C)$. Define $f:=\sum_{i=1}^5 g_i f_i $. Then direct calculation shows that $f$ vanishes at least twice at the 6 points $p_0,\cdots,p_5$ and vanishes on the  15 lines $\overline{p_ip_j}$. As a result, we have proved $h^0(X, \oO_X(-K_X))=1$ for $p_0$ general.

We choose $p_0=[1:2:3:4]$, then $f=8f_1-2f_2-4f_3+1f_4+2f_5$. Using Macaulay2 \cite{M2}, we checked that for $p_0=[1:2:3:4]$ the surface $S'$ has only nodal singularities at the points $p_0,\cdots,p_5$ and is smooth at all the other points. 
Since it requires $A,B,C$ to satisfy finitely many polynomial equations for $S'$ to have singularities at points other than $p_i$ and have singularities other than nodes, we know for $p_0$ general $S'$ has only nodal singularities at the points $p_0,\cdots,p_5$ and is smooth anywhere else. Therefore, blowing-up at $p_i$ resolves the singularities, so the proper transform $S$ of $S'$ in $X$ is smooth. We have a short exact sequence:
\[0\rightarrow \mathcal{O}_X(-{S})\rightarrow\mathcal{O}_X \rightarrow \mathcal{O}_{{S}} \rightarrow 0\]
on $X$, which induces the long exact sequence:
\[\cdots\ra  H^1(X, \mathcal{O}_X)\ra H^1(S, \mathcal{O}_S)\ra H^2(X, \mathcal{O}_X(-S))\ra \cdots.\]
Since $h^1(X, \mathcal{O}_X)$ is a birational invariant, $h^1(X, \mathcal{O}_X)=h^1(\mathbb{P}^3,\mathcal{O}_{\mathbb{P}^3})=0$. By Serre duality, $H^2(X,\oO_X(-S))=H^1(X,\oO_X(K_X+ S))^{\vee}=H^1(X,\mathcal{O}_X)^{\vee}=0$.
One deduces $H^1({S},\mathcal{O}_{{S}})=0$.
Since $S$ is smooth, by the adjunction formula $K_{{S}}=(K_X+{S})\vert_{S}=0$. Hence ${S}$ is a K3 surface.\qed

In the following we assume the six points $p_i$ are general so $S$ is the unique anticanonical section of $X$. 
We observe that there are many rational curves on $S$.
\begin{itemize}
	\item $E_i$ is the exceptional divisor of $S$ over $p_i$, for $i=0,\cdots, 5$. We abuse the notation here and use $E_i$ to represent the exceptional divisors on both $X$ and $S$, where it will be clear from the context whether they are in $X$ or $S$.
	\item $T_{ij}:=E_{ij}\cap S$ is the line $\overline{p_i p_j}$.
	\item There is a unique rational normal curve $R_0$ in $\PP^3$ through the six points $p_0,\cdots, p_5$. Denote by $R$ the proper transform of $R_0$ in $X$.
	\item Let $\Gamma_I$ be the plane in $\PP^3$ through the three points $p_i$ with $i\in I$, for each $I\subset \{0,1,2,3,4,5\}$ such that $\vv{I}=3$.
	\item Set $J= \{0,1,2,3,4,5\}\backslash I$. Then let $L_I$ be the proper transform of the line $\Gamma_I\cap \Gamma_J$ in $X$. By symmetry, $L_I=L_J$, so there are $10$ such lines $L_I$ in $X$. 
\end{itemize}
Then $E_i$ and $T_{ij}$ are all $(-2)$-curves on $S$. We obverse that by direct calculation:
\begin{lemma}
\begin{enumerate}
	\item Each $L_I$ and $R$ are contained in $S$. 
	\item The rational curve $R$ does not meet the ten lines $L_I$.
\end{enumerate}
	\label{fourthline}
\end{lemma}

Let $H_S$ be the restriction of the hyperplane class $H$ to $S$. It follows from Lemma \ref{fourthline} that we have the following relations in $\Pic(S)\cong \NS(S)$: 
\begin{align}
	&H_S\sim E_i+E_j+E_k+T_{ij}+T_{jk}+T_{ik}+L_{ijk}, \quad \txt{for distinct } i,j,k.
	\label{rel4}
\end{align}
Indeed, considering the intersection of $\Gamma_{ijk}$ with $S'$ we have $H_S\sim a(E_i+E_j+E_k)+T_{ij}+T_{jk}+T_{ik}+L_{ijk}$. Since the degree of $S$ is $4$, $H_S^2=4$. Hence $a=1$ by calculating the self-intersection of $H_S$.
By (\ref{rel4}), we can gather the following intersection products over $S$:
\begin{align}
	\begin{split}
	&H_S^2=4, \quad E_i^2=T_{ij}^2=L_{ijk}^2=R^2=-2;\\
	&H_S\cdot E_i=0, \quad H_S\cdot T_{ij}=1, \quad H_S\cdot R=3,\quad H_S\cdot L_{ijk}=1;\\
	&T_{ij}\cdot E_i=1, \quad T_{ij}\cdot L_{ijk}= T_{pq}\cdot L_{ijk}=1, \quad R\cdot E_i=1,
	\label{latticeS}
	\end{split}
\end{align}
for $i,j,k, p, q$ distinct, and all the other intersections among $H_S,E_i,T_{ij},L_{ijk}$ and $R$ are zero.
These intersection products above imply that $\{H_S, E_i, L_I\}$ span a rank-$17$ sublattice of $\NS(S)$. 
\begin{prop}
	The $16$ smooth rational curves $E_i$, $i=0,\cdots, 5$, and $L_I$, $\vv{I}=3$ are pairwise disjoint. Hence $S$ is a K3 Kummer surface.
	\label{166}
\end{prop}
\pf. The only nontrivial part is to show $L_I$ and $L_{I'}$ are disjoint, which by symmetry can be reduced to that $L_{123}$ and $L_{124}$ do not meet. Indeed, if $\{0,1,2,3,4,5\}=I\sqcup J=I'\sqcup J'$ with $\vv{I}=\vv{I'}=3$, then either $\vv{I\cap I'}=2$ or $\vv{I\cap J'}=2$. Now in $\PP^3$, we have $L_{123}\cap L_{124}$ is over $(\Gamma_{123}\cap \Gamma_{045})\cap (\Gamma_{124}\cap \Gamma_{035})$, which equals $(\Gamma_{123}\cap \Gamma_{124})\cap (\Gamma_{045}\cap \Gamma_{035})=\overline{p_1 p_2}\cap \overline{p_0 p_5}=\emptyset$ since $\overline{p_1 p_2}$ and $\overline{p_0 p_5}$ are skew lines. Hence $L_{123}\cap L_{124}=\emptyset$.
Finally, Nikulin's result (Lemma \ref{NK16}) implies that $S$ is a K3 Kummer surface. \qed

By  Lemma \ref{NK16}, the K3 Kummer surface $S$ is associated with an abelian surface $A$ such that $\kum(A)$ is a singular Kummer surface, and there is a natural de-singularization $\pi: S\ra \kum(A)$ at $16$ nodes, such that the $16$ exceptional divisors are exactly those $E_i$ and $L_I$.

\subsection{Generic \texorpdfstring{$S$}{S} has Picard rank 17}
Here we prove that $\rho(S)=17$ when the six points $p_0,\cdots, p_5$ are in very general position.

We recall that an ample line bundle $\mathcal{L}$ on an abelian variety $A$  defines a {\em polarization} $\phi_\mathcal{L}$ which is an isogeny $\phi_\mathcal{L}: A\ra \Pic^0(A)$ sending $x$ to $T_x^* \mathcal{L} \otimes \mathcal{L}\ik$, with $T_x:A\ra A, y\mt y+x$ the translation morphism of adding $x$. The polarization $\varphi_\mathcal{L}$ is principle if it has degree $1$, that is, is an isomorphism.
In the following let $A$ be an abelian surface. 

\begin{lemma}
	Suppose $A$ is an abelian surface with an irreducible effective Cartier divisor $D$ such that $D^2=2$. Then $A$ is principally polarized by $\mathcal{L}=\oO_A(D)$, and $A\cong J(D)$ is the Jacobian variety of $D$.
	\label{Ajacobian}
\end{lemma}
\pf.  Any effective divisor on an abelian variety is nef, so $D$ is nef. Since $D^2=2>0$, $D$ is ample. Now, $\varphi_{\mathcal{L}}$ is a polarization on $A$. Since $\deg \varphi_{\mathcal{L}}= D^2/2=1$ (see \cite[Thm. 5.2.4]{Birkenhake2004}), we have $D$ gives a principal polarization.
Finally, since $D$ is irreducible, $A\cong J(D)$ by the Matsusaka-Ran Criterion (see \cite[11.8.1]{Birkenhake2004}).\qed

Now we fix $S$ to be the unique anticanonical section of $X$. Let $A$ be the abelian surface which $S$ is associated with. Consider the commutative diagram (\ref{cd}), where $\tilde{S}$ is the double cover of $S$. We let $C$ be the double cover of $R$ in $\tilde{S}$, and let $C_0$ be the image of $C$ in the abelian surface $A$. Recall that $S$ is of Jacobian type if $A$ is the Jacobian of a smooth genus-$2$ curve.
\begin{prop}
	The abelian surface $A$ is isomorphic to $J(C_0)$, so $S$ is of Jacobian type.\label{Sjacobian}
\end{prop}
\pf. By Lemma \ref{fourthline} (2), $R$ does not meet the $10$ rational curves $L_{ijk}$. Hence the double cover $C\ra R$ ramifies at six distinct double points corresponding to $p_i$. Hence $C$ is a smooth genus-$2$ curve. Now $\pi':\tilde{S}\ra A$ is the blow-up of $A$ at $16$ smooth points, and $R\cdot E_i=1$. Hence $C_0\cong C$ is a smooth genus $2$ curve. By Lemma \ref{Ajacobian}, we need only show $C_0^2=2$. Indeed, the $16$ exceptional divisors of the smooth blow-up $\pi'$ are just $E_i$ and $L_I$ for $0\leq i\leq 5$ and $\vv{I}=3$, which are the branch loci of the double-cover $f':\tilde{S}\ra S$. Then we have $\pi'^*C_0\sim C+\sum_{i=0}^5 E_i$. Hence  $C_0\cdot C_0=(\pi'^*C_0)\cdot(\pi'^*C_0)=C^2+2C\cdot\sum_{i=0}^5 E_i+\sum_{i=0}^5 (-1)=C^2+2\cdot 6-6=C^2+6$. So we only need to show that $C^2=-4$. Consider again the double cover $f'$. We have $K_{\tilde{S}}\sim f'^* K_S+\sum_{i=0}^5 E_i+\sum_{\vv{I}=3} L_I=\sum_{i=0}^5 E_i+\sum_{\vv{I}=3} L_I$. By Lemma \ref{fourthline}, $C$ does not meet $L_I$, and $C\cdot E_i=1$ in $\tilde{S}$. Hence $K_{\tilde{S}}\cdot C=(\sum_{i} E_i+\sum_I L_I)\cdot C=6$. By the adjunction formula on $\tilde{S}$ we have $(K_{\tilde{S}}+C)\cdot C=\deg K_C=2$.  Hence $C^2=2-6=-4$. This shows that $C_0^2=2$, so $A=J(C_0)$. \qed

\begin{prop}
	For the six points $p_0,\cdots, p_5$ in very general position, $\rho(S)=17$.\label{S17}
\end{prop}
\pf. By Formula (\ref{add16}) and Proposition \ref{Sjacobian},  we need only show $\rho(J(C_0))=1$ for the six points in very general position. Let $M_2$ be the moduli space of smooth genus $2$ curves. By \cite[IV. Ex. 2.2]{Hartshorne1977}, there is an isomorphism $P:M_{2}\ra M_{0,6}/\mathcal{S}_6$ sending a smooth genus $2$ curve $C$ to the six points on $\PP^1$ over which $C\ra \PP^1$ branches. Now the moduli space $M_{0,6}$ is naturally isomorphic to the moduli of rational normal curves in $\PP^3$, which is in turn the moduli of six points in $\PP^3$ in linearly general position. Hence if the six points $p_0,\cdots, p_5$ in very general position, then the corresponding double cover $C_0$ is a very general genus $2$ curve, so that $\rho (J(C_0))=1$ by \cite{koizumi1976}.\qed

\subsection{The Jacobian Kummer structure} In this paragraph, we assume the six points are in general position, so that $S$ is a Jacobian K3 Kummer surface associated with $A=J(C_0)$. Define $\mathcal{N}_{EL}:=\{E_i\mid i=0,\cdots, 5\} \cup \{L_I\mid \vv{I}=3\}$. 

We follow \cite[1.5, 1.6]{Keum1997}. The embedding of $C_0$ into $A$ realizes $C_0$ as a theta divisor $\Theta$. There are exactly six order-$2$ points $x_0, \cdots, x_5$ on $C_0$, which correspond to the six points $p_i$ in $\PP^3$. Therefore, if we fix a choice of the identity among $x_i$, say $x_0$, then the $16$ points in $A[2]$ are:
\begin{align*}
	\mu_i=[x_i-x_0],\quad 0\leq i\leq 5, \quad \mu_{jk}=[x_j+x_k-2x_0], \quad 1\leq j<k\leq 5.
\end{align*}
Moreover, we identify $\mu_i$ with $i$ and $\mu_{jk}$ with $jk$. Then
\[A[2]=\{i,jk\mid 0\leq i\leq 5, 1\leq j<k\leq 5\}.\]
Under this identification, the group law on $A[2]$ is given by
\begin{align}\label{law}
	2i=2jk=0,\; i+j=ij,\; i+0=i,\; jk+0=jk,\; jk+jm=km,\; jk+mp=q.
\end{align}
where $\{1,2,3,4,5\}=\{j,k,m,p,q\}$.

Now it is clear that we can identify $E_i$ with $i$ by permuting the six points $p_i$. The question is to correctly identify those $L_I$ with $jk$. We have:
\begin{theorem}\label{kumaddition}
If we choose $x_0$ to be the identity and identify $A[2]$ with $\{i,jk\}$ as above and identify each $E_i$ with $i$. Then $L_{0jk}=jk$ for $1\leq j<k\leq 5$. That is, the bijection between $\mathcal{N}_{EL}$ with $A[2]$ is given by
\[ E_i\mt i,\quad L_{0jk}\mt jk.\]
\end{theorem}

We prove Theorem \ref{kumaddition} by finding the hyperplanes in $A[2]$.
\begin{lemma}	\label{KS}
	Suppose $\{0,1,2,3,4,5\}=\{i,j,k,m,p,q\}$. Then the following $30$ classes are in $\NS(S)$:
	\begin{align}	
		C_{ij}:=&\frac{1}{2}\left(E_k+E_m+E_p+E_q+L_{ijk}+L_{ijm}+L_{ijp}+L_{ijq}\right),\\
		\label{CCij}
		D_{ij}:=&\frac{1}{2}\sum_{\xi\in\mathcal{N}_{EL}}\xi-C_{ij}.
	\end{align}
\end{lemma}
\pf.
By Formula (\ref{rel4}),
\begin{align*}
	&L_{ijk}+L_{ijm}+L_{ijp}+L_{ijq}=L_{mpq}+L_{kpq}+L_{kmq}+L_{kmp}\\
	=&4H-3(E_k+E_m+E_p+E_q)-2(T_{km}+T_{kp}+T_{kq}+T_{mp}+T_{mq}+T_{pq}).
\end{align*}
As a result,
\[2C_{ij}= 4H-2(E_k+E_m+E_p+E_q+T_{km}+T_{kp}+T_{kq}+T_{mp}+T_{mq}+T_{pq}).\]
Hence $C_{ij}\in \NS(S)$.  Next, we can directly compute that $(1/2)\sum_{\xi\in\mathcal{N}_{EL}}\xi\in \NS(S)$ by Formula (\ref{rel4}), or recall Nikulin's result that if $S$ is a complex K3 surface with $16$ disjoint smooth rational curves $C_i$, then $(1/2)\sum_i C_i\in \NS(S)$. (see \cite{Nikulin1975}. Also see \cite[\S 14, Rem 3.16]{Huybrechts2016}) Hence $D_{ij}=(1/2)\sum_{\xi\in\mathcal{N}_{EL}}\xi-C_{ij}\in \NS(S)$.\qed

\pfof{Theorem \ref{kumaddition}}.
By Lemma \ref{KS} and Lemma \ref{NK2}, we know the $30$ hyperplanes of $8$-elements in $\mathcal{N}_{EL}$ are the following:
\begin{align*}
	\Gamma_{ij}&=\{E_i,E_j,L_{ikm},L_{ikp},L_{ikq},L_{imp},L_{imq},L_{ipq}\},\\
	\Gamma_{ij}^c&=\{E_k,E_m,E_p,E_q,L_{ijk},L_{ijm},L_{ijp},L_{ijq}\},
\end{align*}
for $0\leq i<j \leq 5$. Therefore, intersecting every pair of hyperplanes gives us all the $140$ affine $2$-planes in $A[2]\cong \F_2^{\oplus 4}$. They are:
\[  \begin{array}{ll}
\{E_i,E_j,E_k,L_{ijk}\}, &\{E_i,E_j,L_{ikm},L_{jkm}\},\\
	\{E_i,L_{ijk},L_{ijm},L_{ikm}\}, &\{L_{ikp},L_{ikq},L_{imp},L_{imq}\},
 \end{array} \] 
for $i,j,k,m,p,q\in \{0,\cdots, 5\}$ distinct.
Note that the description of the affine $2$-planes does not rely on the identity point we choose.
Recall a simple fact that in an affine $2$-plane $\{r_1,r_2,r_3,r_4\}$ over $\F_2$, if one of $r_i=0$, say $r_1$, then $r_2+r_3=r_4$. Since we choose $E_0\mt 0$ to be the identity, we can apply this fact to all the $2$-planes above, and write the addition in $\mathcal{N}_{EL}$ by $\oplus$. Then we find:
\begin{enumerate}
	\item The $2$-plane $\{E_0, E_j,E_k, L_{0jk}\}$ gives $L_{0jk}=E_j\oplus E_k$.
	\item The $2$-plane $\{E_0, E_j,L_{0km}, L_{jkm}\}$ gives $E_j=L_{0km}\oplus L_{jkm}=L_{0km}\oplus L_{0pq}$.
	\item The $2$-plane $\{E_0,L_{0jk},L_{0jm}, L_{0km}\}$ gives $L_{0jk}=L_{0jm}\oplus L_{0km}$.
\end{enumerate}
By our discussion above, when fixing $E_0\mt 0$,  the group structure on $\mathcal{N}_{EL}$ is exactly given by (\ref{law}). Since we already identified $E_i$ with $i$, we must have $L_{0jk}\mt jk$.\qed

\begin{remark}
	This proof does not assume $\rho(S)=17$.	\label{17}
\end{remark}

\section{Keum's 192 Automorphisms}\label{k192}
In the rest paragraphs of the paper, we assume the six points $p_i$ are in very general position, so that $\rho(S)=17$. Let $A$ be the abelian surface associated with $S$. 
We will always identify $\mathcal{N}_{EL}$ with $A[2]$ via Theorem \ref{kumaddition}, that is, $E_i=N_i$ and $L_{0jk}=N_{jk}$. 
\subsection{The hyperplane section $\Lambda$}
When $A=J(C)$ is the Jacobian of a smooth genus $2$ curve $C$, $\kum(A)$ embeds in $\PP^3$ as a quartic surface with exactly $16$ nodes (see \cite[3.1]{Keum1997}). Let $\Lambda'$ be the hyperplane class of $\kum(A)$ under such an embedding, and let $\Lambda=\pi^*\Lambda'$ on $S$. 

We note that Keum wrote $H$ for our $\Lambda$. Also, $\Lambda\not\equiv H_S$ where we denote by $H_S$ the hyperplane class of $S$ from $S'$. In particular, $S'\not\cong \kum(A)$ since the former has only $6$ nodes, while $\kum(A)$ has $16$.

Since $\rho(S)=17$, we have $\Lambda^2=4$, $\Lambda\cdot E_i=\Lambda\cdot L_{ijk}=0$, and $\{\Lambda, E_i, L_I\}$ freely generate $\NS(S)_\Q$ (see \cite[Thm. 1]{Naruki1991}). 
On the other hand, by the intersection products in (\ref{latticeS}), we know that $H_S\not\in \Z\mathcal{N}_{EL}$. Hence $\{H_S, E_i, L_I\}$ also generate $\NS(S)_\Q$. Now we describe the relations between these classes:
\begin{prop}
	\label{LH}
In $\NS(S)$, we have
\begin{align*}
H_S&\sim \frac{1}{2}\left(3\Lambda-\sum_{\vv{I}=3} L_I\right),\\
T_{ij}&\sim \frac{1}{2}(\Lambda-E_i-E_j-\sum_{k\neq i,j} L_{ijk}),\\
R&\sim \frac{1}{2}(\Lambda-\sum_i E_i).
\end{align*}
Furthermore, $\Lambda\cdot H_S=6$, $\Lambda\cdot T_{ij}=2$, and $\{H_S,E_i, T_{ij}\}$ also generate $\NS(S)_\Q$. 
\end{prop}
\pf. Since $S$ is K3 and $\rho(S)=17$, The first equation we want to prove is equivalent to $3\Lambda\sim 2H_S+\sum_{\vv{I}=3} L_I$.
We claim that the only $\Q$-divisor $D\in \Pic(S)_\Q$ with $D\cdot E_i=0$ and $D\cdot L_{ijk}=0$ is of the form $D\sim  r\left(2H_S+\sum_{\vv{I}=3} L_I \right)$ for some $r\in \Q$. Indeed, we have shown above that $\{H_S, E_i, L_I\}$ also generate $\NS(S)_\Q$.  Suppose $D\sim hH_S-\sum_i a_i E_i-\sum b_{ijk}L_{ijk}$. By Formulas (\ref{latticeS}), we have $0=D\cdot E_i=2a_i$ and $0=D\cdot L_{ijk}=h+2b_{ijk}$, so $a_i=0$ and $b_{ijk}=-h/2$. Hence $D\sim (h/2)\left(2H_S+\sum_I L_I\right)$.
Next, $D^2=r^2(2H_S+\sum_I L_I)^2=36r^2$. Therefore requiring $D^2=4$ gives $r=\pm 1/3$. Clearly, $D$ is effective if and only if $r=1/3$.

By checking numerical equivalence we can prove the equalities for $T_{ij}$ and $R$. Since  $\Lambda^2=4$ and $\Lambda\cdot E_i=\Lambda\cdot L_{ijk}=0$, we have $\Lambda\cdot H_S=\Lambda\cdot (3/2)\Lambda=6$, and $\Lambda\cdot T_{ij}=\Lambda^2/2=2$.
Finally, to show $\{H_S,E_i, T_{ij}\}$ also generate $\NS(S)_\Q$, we need only show every $L_I$ is generated by $H_S,E_i$ and $T_{ij}$ over $\Q$, which follows from Formula (\ref{rel4}). 
\qed

\begin{remark}
	A K3 Kummer surface associated with $A$ has $16$ $(-2)$-curves called tropes. The $16$ tropes and the $16$ nodes $N_\alpha$ form the $(16,6)$-configuration where every trope passes through six nodes, and every node is on six tropes. Our $T_{ij}$ and $R$ are exactly the $16$ tropes. Our notations translate to Keum's \cite[1.8]{Keum1997}  (also see \cite{Naruki1991}) as follows: Our $R$ to $T_0$, $T_{0j}$ to $T_j$ and $T_{jk}$ to $T_{jk}$ for $j,k=1,2,3,4,5$.
	\label{tropes}
\end{remark}

\subsection{Keum's $192$ automorphisms}
By \cite[Def. 6.12]{Keum1997}, a Weber Hexad $\mathcal{H}$, as a subset of $A[2]$, is the symmetric difference of a G\"{o}pel tetrad with a Rosenhain tetrad (see \cite[\S 2]{Keum1997} and \cite[10.2]{Birkenhake2004}). There are $192$ Weber Hexads. The translations $t_\alpha: A[2]\ra A[2]$, $x\mt x+\alpha$ send a Weber Hexad to another Weber Hexad. In particular, one of the Weber Hexads is
\[\mathcal{H}_1=\{0,14,15,23,25,34\}.\]
\begin{theorem} \cite[Thm. 6.11, 6.16]{Keum1997}
	For any Weber Hexad $\mathcal{H}$, the complete linear system 
	\[\left|7\Lambda-4\sum_{h\in \mathcal{H}} N_h \right|\]
	induces an automorphism $\kappa_\mathcal{H}:S\ra S$ of infinite order. 
	\label{KeumH}
	The automorphism $\kappa=\kappa_\mathcal{H}$ is determined by its action on the Picard lattice: $\kappa^*: \NS(S)\ra \NS(S)$. For instance the Hexad $\mathcal{H}_1=\{0,14,15,23,25,34\}$ gives the automorphism $\kappa_1$ such that
	\begin{gather*}
	\kappa_1^*  \Lambda=7\Lambda -4(N_0+N_{14}+N_{15}+N_{23}+N_{25}+N_{34}),  \quad \kappa_1^* N_{12}=N_{12},\\
	\kappa_1^* N_2=N_3, \quad \kappa_1^* N_3=N_{13}, \quad 	\kappa_1^* N_{13}=N_2, \\
	\kappa_1^* N_1=N_4, \quad \kappa_1^* N_4=N_{24}, \quad 	\kappa_1^* N_{24}=N_1, \\
	\kappa_1^* N_{14}=M-N_{25}, \quad \kappa_1^* N_{23}=M-N_{15}, \quad \kappa_1^* N_{5}=M-N_{34}, \\
	\kappa_1^* N_{34}=M-N_{0}, \quad \kappa_1^* N_{35}=M-N_{23}, \quad \kappa_1^* N_{45}=M-N_{14}, \\
	\kappa_1^* N_0=N_5, \quad \kappa_1^* N_{15}=N_{35}, \quad \kappa_1^* N_{25}=N_{45},
\end{gather*}
where $M:=2\Lambda -\sum_{h\in \mathcal{H}_1}N_h=2\Lambda- (N_0+N_{14}+N_{15}+N_{23}+N_{25}+N_{34})$. 
\end{theorem}

One of the goals of the paper is to find a pseudo-automorphism on $X$ which restricts to Keum's automorphisms. Instead of $\mathcal{H}_1$ above, we consider a different Weber Hexad. The translation map on $A$ by the point $x_5$ induces an automorphism (also called the translation) $t_5$ of $S$ \cite[5.1(i)]{Keum1997}. Let $\mathcal{H}:=t_5(\mathcal{H}_1)=\{5,23,1,14,2,12\}$. Then $\mathcal{H}$ is a Weber Hexad. Then $t_5\circ \kappa_1 \circ t_5=\kappa_\mathcal{H}$ (To see this, we can show their actions on $\NS(S)$ agree and then use Proposition \ref{autS}). 

In the rest of the paper, we will let $\kappa$ be the automorphism associated with $\mathcal{H}=\{5,23,1,14,2,12\}$. Then $\kappa^*$ is given by:
\begin{align}
	\begin{split}
	\begin{gathered}
\kappa^*  \Lambda=7\Lambda -4(N_1+N_{2}+N_{5}+N_{12}+N_{23}+N_{14}),  \quad \kappa^* N_{34}=N_{34},\\
\kappa^* N_{25}=N_{35}, \quad \kappa^* N_{35}=N_{24}, \quad 	\kappa^* N_{24}=N_{25}, \\
\kappa^* N_{15}=N_{45}, \quad \kappa^* N_{45}=N_{13}, \quad 	\kappa^* N_{13}=N_{15}, \\
	\kappa^* N_{23}=U-N_{2}, \quad \kappa^* N_{14}=U-N_{1}, \quad \kappa^* N_{0}=U-N_{12}, \\
	\kappa^* N_{12}=U-N_{5}, \quad \kappa^* N_{3}=U-N_{14}, \quad \kappa^* N_{4}=U-N_{23}, \\
	\kappa^* N_5=N_0, \quad \kappa^* N_{1}=N_{3}, \quad \kappa^* N_{2}=N_{4},
\end{gathered}
	\end{split}
	\label{mu0}
\end{align}
where $U:= 2\Lambda -\sum_{h\in \mathcal{H}}N_h=2\Lambda - (N_1+N_{2}+N_{5}+N_{12}+N_{23}+N_{14})$.

\begin{remark}
	Kond\={o} \cite{Kondo1998} proved that $\Aut(S)$ of a general Jacobian K3 Kummer $S$ is generated by the `classical' automorphisms along with Keum's 192 automorphisms of infinite order.	\label{kondo}
\end{remark}

\begin{corollary}
	\label{kappaTR}
	In $\NS(S)$, we have $\kappa^*R=R$. 
\end{corollary}
\pf. This follows from a calculation using Formula  (\ref{mu0}) and Proposition \ref{LH}.\qed

We conclude this section by two auxiliary results.
First, consider the restriction map $r:\Pic(X)\ra \Pic(S)$. Then $r(H)=H_S$, $r(E_i)=E_i$, and $r(E_{ij})=T_{ij}$. Recall that if  $D\sim dH-\sum_i m_i E_i -\sum_{ij} m_{ij} E_{ij}\in \Pic(X)$, then we say $\deg D=d$. We have
\begin{lemma}
	\label{rkernel}
	Suppose $D_1, D_2\in \Pic(X)$ such that $r(D_1)\sim r(D_2)$. Then $\deg D_1=\deg D_2$.
\end{lemma}
\pf. Since $r$ is linear, we need only show that if $r(D_1)\sim 0$, then $\deg D_1=0$. Suppose $D_1\sim dH-\sum_i m_i E_i -\sum_{ij} m_{ij} E_{ij}$ and $r(D_1)\sim 0$. Then $r(D_1)\cdot \Lambda=r(D_1)\cdot L_I=0$ for each $i$ and $ I$. By (\ref{latticeS}) and Proposition \ref{LH}, $0=r(D_1)\cdot \Lambda=6d-2\sum_{ij} m_{ij}$, and $
r(D_1)\cdot L_{ijk}=d-(m_{ij}+m_{ik}+m_{jk}+m_{pt}+m_{pq}+m_{qt})$ for $\{0,1,2,3,4,5\}=\{i,j,k,p,q,t\}$. Therefore $0=r(D_1)\cdot \sum_I L_I=10d-4\sum_{ij} m_{ij}$. This makes $d=0$.\qed

\begin{prop}
	\label{nefbigH}
	Both $\Lambda$ and $H_S$ are nef and big, with $h^0(S,\oO_S(\Lambda))=h^0(S,\oO_S(H_S))=4$.
\end{prop}
\pf. First $\Lambda$ is the pullback of the hyperplane class via the embedding $\kum(A)\hookrightarrow \PP^3$. Then $H_S$ is the pullback of $H$ via the embedding $S\hookrightarrow X$. Since both hyperplane classes here are ample, $\Lambda$ and $H_S$ are big and nef. Since $S$ is K3, by the Kawamata-Viehweg vanishing theorem,  $H^i(S, \oO_S(H_S))=0$ for $i>0$. Hence the Riemann-Roch theorem implies that $h^0(S, \oO_S(H_S))=(1/2)H_S^2+2=4$. A similar reasoning shows $h^0(S,\oO_S(\Lambda))=(1/2)\Lambda^2+2=4$. \qed

\section{The 9 Rational Quartics}\label{9Q}
In this section we define $9$ quartics in $X$ and $Y$ such that they restrict to some of the $(-2)$-curves appearing in the mapping table (\ref{mu0}) of Keum's automorphism $\kappa$.
\begin{defi}\label{DQ}
Define
\begin{align*}	D:=\quad &  13H-7(E_1+E_2+E_5)-5(E_0+E_3+E_4)\\
	&-3(E_{03}+E_{04}+E_{34})
	-4(E_{05}+E_{13}+E_{24})-(E_{12}+E_{15}+E_{25}).
\end{align*}
Define the following $9$ quartic classes in $Y$:
\[	\begin{array}[h]{lll}
	Q_0 :=& 4H-2E_0-E_3-E_4-2E_1-2E_2-3E_5\\
	&\hspace{5em}-E_{03}-E_{04}-2E_{05}-E_{13}-E_{24} -E_{15}-E_{25},\\
	Q_3 :=& 4H-E_0-2E_3-E_4-3E_1-2E_2-2E_5\\
	&\hspace{5em}-E_{03}-E_{34}-E_{05}-2E_{13}-E_{24}-E_{12}-E_{15}, \\
	Q_4 :=& 4H-E_0-E_3-2E_4-2E_1-3E_2-2E_5\\
	&\hspace{5em}-E_{04}-E_{34}-E_{05}-E_{13}-2E_{24} -E_{12}-E_{25}, \\
	Q_{05} :=& A-E_{34},\qquad Q_{13} := A-E_{04}, \qquad Q_{24} := A-E_{03},\\
	Q_{12} :=& A-E_{05},\qquad Q_{15} := A-E_{24}, \qquad Q_{25} := A-E_{13},
	\end{array}
\]
where $A:=4H-2\sum_{i=0}^5 E_i-(E_{03}+E_{04}+E_{34})-(E_{05}+E_{13}+E_{24})$. We will refer to these quartics as $Q_\alpha$, where \[\alpha\in \mathcal{A}:=\{0,3,4,12,15,25,05,13,24\}.\] Finally, define $Q_\alpha'$ to be the image of $Q_\alpha$ in $\PP^3$.
\end{defi}

\begin{remark}
\begin{enumerate}
	\item Consider the action of $\mathcal{S}_3$ on $\Pic(X)$ by permuting the ordered pairs of points $\{(p_5,p_0),(p_1,p_3),(p_2,p_4)\}$. Then $D$ and $A$ are fixed by this action. We can divide  $\{Q_\alpha\}$ into three subsets: $\{0,3,4\}\cup \{05,13,24\}\cup\{12,15,25\}$ where each subset is contained in an orbit of the $\mathcal{S}_3$-action. 
	\item Since each $Q_\alpha$ is effective (Theorem \ref{pencilsQ}), each $Q_\alpha '$ is a singular quartic in $\PP^3$ with either a triple point (for $\alpha=0,3,4$) or a double line (for all of the nine). A singular quartic $Q$ in $\PP^3$ with a triple point or a double line must be rational \cite{Jessop1916}. Therefore all the $Q_\alpha$ are rational. In Theorem \ref{psdautoY} we show $\phi$ is a pseudo-automorphism. Then we find $\phi:Q_\alpha\dashrightarrow E_\alpha$ is birational, which shows that those $Q_\alpha$ are rational in a different way. Each $Q_\alpha$ spans an extremal ray of $\Effb(Y)$ and $\Effb(X)$ since $Q_\alpha=\phi^* E_\alpha$ and $E_\alpha$ is extremal.
\end{enumerate}
	\label{qrmk}
\end{remark}
 
\begin{theorem}\label{pencilsQ}
	For the six points $p_0,\cdots, p_5$ in very general position, we have
	\begin{enumerate}
		\item Consider the restriction map $r:\Pic(X)\ra \Pic(S)$. Then $r(D)=\kappa^* H_S$, and
	\[\begin{array}{ll}
		& r(Q_0)=\kappa^* N_0=U-N_{12},\quad  r(Q_3)=\kappa^* N_3=U-N_{14},\quad  r(Q_4)=\kappa^* N_4=U-N_{23},\\
		& r(Q_\alpha)=\kappa^* T_{\alpha} \txt{ for }\alpha\in \{05,13,24,12,15,25\}.
		\end{array}\]
		where $U:= \Lambda - (N_1+N_{2}+N_{5}+N_{12}+N_{23}+N_{14})$ as in (\ref{mu0}).
	\item  Consider $Q_\alpha$ as divisor classes on $X$. For each $\alpha\in \mathcal{A}$, $h^0(X,\oO_X(Q_\alpha))=1$. The unique global sections of $Q_\alpha$ are irreducible and distinct from each other.
	\item
		$h^0(X,\oO_X(Q_\alpha-E_i))=0$, and $h^0(X,\oO_X(Q_\alpha-E_{ij}))=0$, for each $0\leq i\neq j\leq 5$.
	\end{enumerate}
\end{theorem}

\pf. (1) When the six points $p_i$ are very general, $\rho(S)=17$, and $\{\Lambda,E_i, L_{0jk}\}$ generate $\Pic(S)_\Q$.  Since $S$ is K3, for each equality we need only show that the intersection products of both sides with the $\Q$-basis $\{\Lambda,E_i, L_{0jk}\}$ coincide. 
We use that $r(H)=H_S$, $r(E_i)=E_i$, $r(E_{ij})=T_{ij}$, the intersection products from  (\ref{rel4}), (\ref{latticeS}) and Proposition \ref{LH}. 
Then (1) follows from a direct calculation.

(2) For each $Q_\alpha$, there exists an exact sequence:
\[0\ra H^0(X,\oO_X(Q_\alpha-S))\ra H^0(X,\oO_X(Q_\alpha))\ra H^0(S, \oO_S({Q_\alpha} {\mid_S})).\]
By definition, it is easy to verify that each $Q_\alpha-S\not\sim 0$, has degree zero, but with negative coefficients on some $E_{ij}$. Hence none of those $Q_\alpha-S$ are effective. Thus $h^0(X,\oO_X(Q_\alpha-S))=0$. On the other hand, by (1) we know each $Q_\alpha$ restricts to the preimage of a $(-2)$-curve on $S$ under $\kappa$, which is also a $(-2)$-curve, so that $h^0(S, \oO_S({Q_\alpha} {\mid_S})=1$. Therefore, $h^0(X,\oO_X(Q_\alpha))\leq 1$ for each $\alpha\in \mathcal{A}$.

	It remains to show that each $Q_\alpha$ is indeed effective.
	Here we let $[x:y:z:w]$ be the homogeneous coordinates on $\PP^3$ and let five of the six points be at standard position and the sixth at $[1:a:b:c]$ for general $a,b,c$. To make the polynomials simpler we will choose different orders of the six points for each case. By symmetry, we need only show $Q_0, Q_{24}$ and $Q_{12}$ (we choose $Q_{24}$ and $Q_{12}$ also for the sake of the proof of Theorem \ref{E4image}). We claim the polynomials $f_0,f_{24},f_{12}$ defining $Q'_0, Q'_{24}$ and $Q'_{12}$ are:

	(i) For $Q_0$, we place $(p_0,\cdots, p_5)$ at 
	$([1: 0: 0: 0], [0: 1: 0: 0], [0: 0: 1: 0], [1: a: b: c], [1: 1: 1: 1], [0: 0: 0: 1])$. 
	Then it is easy to see the following polynomials define sections of $4H-2E_0-2E_1-2E_2-E_4-3E_5-E_{04}-2E_{05}-E_{24}-E_{15}-E_{25}$:
	\[(x y z (y-w),y z^2 (y-w),x z^2 (y-w),x y (x-w) (y-z),x z (x-y) (y-z),y z (z-w)  (x-y)).\]
	Let 
	\begin{align}\label{f0}
	\begin{split}
		f_0&:= b c (-a+b-1)x y z (y-w)+a (c-b)y z^2 (y-w)+a bx z^2 (y-w)\\
	&+ b^2 c x y (x-w) (y-z)-a b c x z (x-y) (y-z)+b (c-a)y z (z-w)  (x-y).
	\end{split}
	\end{align}
	Then it is easy to verify that $f_0$ vanishes at $p_3$, line $03$ and line $13$. Hence $f_0$ defines the unique quartic $Q_0$.

	(ii) For $Q_{12}$ and $Q_{24}$, we place $(p_0,\cdots, p_5)$ at 
	$([0: 1: 0: 0], [1: 0: 0: 0], [1: a: b: c], [0: 0: 1: 0], [0: 0: 0: 1], [1: 1: 1: 1])$.
Then the polynomial
\begin{align}\label{f12}
	\begin{split}
	f_{12}&:=-a ( b-1)^2 c yz (x-w )  (z-w )+ a (b - c) c yz  ( x-w)  (x - z) \\
	&-b (b - c) c y^2 (x-w)  (x - z)+ a b (1 - 2 c +   b c)  w y (x - z) ( z-w)\\
	&+ a^2 b ( c-1) w x (x - z) ( z-w) -a ( b-1) b c x y (x - z) ( z-w).
	\end{split}
\end{align}
defines the unique quartic $Q_{12}$. Similarly, the polynomial
\begin{align}\label{f24}
	\begin{split}
	f_{24}&:=(a - b) (a - c) (b - c)yw (x - z) (x - w) + a (a - b) b (a - c) ( c-1) xw (x - z) (y - w)\\
	&-( a-1) a b (a - c) ( c-1)  xw (x - z) (y - z)+ a (a - b) ( b-1) (c-1) c x z (x - w) (y - w)\\
	&+(a-1) (b-1) b (b - c) c x y (x - w) (y - w) -a (b-1)^2 (a - c) c x y (x - w) (z - w).
 	\end{split}
\end{align}
defines the unique quartic $Q_{24}$.

Now each $Q_\alpha$ is distinct because their restriction to $S$ are distinct $(-2)$-curves. Finally, we show each $Q_\alpha$ is irreducible. We fix $\alpha$ and suppose $Q_\alpha$ is reducible. Then $Q_\alpha\sim D_1+D_2$ for $D_1$ and $D_2$ both nontrivial and effective. Then we find $r(D_1)+r(D_2)\sim r(Q_\alpha)$ is a $(-2)$-curve, hence irreducible and not a sum of two nontrivial effective classes. This implies that either $r(D_1)$ or $r(D_2)$ is trivial or not effective. Suppose $r(D_1)\sim 0$. By Lemma \ref{rkernel}, $\deg D_1=0$. Otherwise, suppose $r(D_1)$ is not effective. Since $S$ is irreducible (Proposition \ref{uniqueminusK}), we conclude that $S$ is contained in the fixed part of $D_1$, which implies that $\deg D_1\geq 4$. Since $\deg Q_\alpha=4$, we must have $\deg D_1=4$ and $\deg D_2=0$. As a conclusion, in either case, one of $D_1$ and $D_2$ must have degree $0$. Assume $\deg D_1=0$. Then $D_1$ is an effective sum of some $E_i$ and $E_{ij}$. Now $Q_\alpha-D_1=D_2$ is effective, which contradicts (3) proved in the following. Hence $Q_\alpha$ is irreducible.

(3). Consider a similar exact sequence:
\[0\ra H^0(X,\oO_X(Q_\alpha-E_i-S))\ra H^0(X,\oO_X(Q_\alpha-E_i))\ra H^0(S, \oO_S((Q_\alpha-E_i) {\mid_S})).\]
Then $h^0(X,\oO_X(Q_\alpha-E_i-S))=0$. So we need only show that $h^0(S, \oO_S((Q_\alpha-E_i){\mid_S}))=0$. Indeed, each ${Q_\alpha}_{\mid_S}$ equals to the preimage of a $(-2)$-curve under $\kappa$, hence is a $(-2)$-curve. By (1) and some calculations, we find ${Q_\alpha}_{\mid_S}\neq E_i$ or $T_{ij}$. Now over the K3 surface $S$, if $C_1$ and $C_2$ are two distinct $(-2)$-curves, then $C_1-C_2$ is not effective. Hence $h^0(S, \oO_S((Q_\alpha-E_i) {\mid_S}))=0$ and $h^0(S, \oO_S((Q_\alpha-E_{ij}) {\mid_S}))=0$. Therefore $H^0(X,\oO_X(Q_\alpha-E_i))=0$, and similarly, $H^0(X,\oO_X(Q_\alpha-E_{ij}))=0$.
\qed

\section{The Linear System \texorpdfstring{$\vv{D}$}{|D|}}\label{dD}
In this section we consider the complete linear system $\vv{D}$. We show $\dim \vv{D}=3$. Then we construct various sections of $\vv{D}$ which arise from planes $\Gamma_{ijk}$ and the $9$ quartics $Q_\alpha$.
\subsection{A first choice}
By Theorem \ref{pencilsQ}, we can make the following definitions:

\begin{defi}Up to nonzero scalars,
\begin{enumerate} 
    \item let $f_\alpha$ be the irreducible quartic polynomial  defining $Q_\alpha'$ in $\PP^3$;
\item let $p_{ijk}$ the linear polynomial defining $\Gamma_{ijk}$, the plane in $\PP^3$ through the points $p_i,p_j$ and $p_k$.
\item Define $4$ polynomials of degree $13$:
\begin{equation}
(s_0,s_1,s_2,s_3):= (p_{034} f_0 f_3 f_4, p_{045} f_3 f_4 f_{24},  p_{234} f_0 f_3 f_{13}, p_{013} f_0 f_4 f_{05}).
		\label{psisection}
	\end{equation}
\item Let $\bar{s}_i$ be the proper transform of the zero locus of $s_i$ in $Y$. Similarly define $\overline{p}_{ijk}$ and $\overline{f}_\alpha$.
\end{enumerate}
\end{defi}
We note that by Theorem \ref{pencilsQ}, $\overline{f}_\alpha$ is the unique section of $Q_\alpha$. Also  $\overline{p}_{ijk}$ is the unique section of $\tilde{\Gamma}_{ijk}$, the proper transform of $\Gamma_{ijk}$.

\begin{prop}\label{wd}
	For the six points in very general position, we have $h^0(X,\oO_X(D))=\linebreak h^0(Y,\oO_Y(D))=4$. The polynomials $s_0,s_1,s_2,s_3$
	are linearly independent. Let $x_{E_i}$ and $x_{E_{ij}}$ be the unique section in $\vv{E_i}$ and $\vv{E_{ij}}$ in $Y$. Then 
	\begin{equation}
		(\bar{s}_0 x_{E_{12}}x_{E_{15}}x_{E_{25}},\bar{s}_1 x_{E_{4}}x_{E_{12}},\bar{s}_2 x_{E_{3}}x_{E_{15}},\bar{s}_3 x_{E_{0}}x_{E_{25}})
		\label{Dsection}
	\end{equation}span the complete linear subsystem $\vv{D}$ over $Y$. 
\end{prop}
\pf. Suppose $s_i$ in (\ref{psisection}) are not linearly independent. Then there are not-all-zero constants $a_i$ such that $a_0s_0+a_1 s_1+a_2 s_2+a_3 s_3=0$. By definition, $f_0\mid s_0,s_2$ and $s_3$, so we must have $f_0\mid a_1 s_1$. By Theorem \ref{pencilsQ} (2), those $Q_\alpha$ are distinct and irreducible. Hence $f_0\nmid s_1$, so $a_1=0$. Repeat for $f_3$ and $f_4$ we have $a_2=a_3=0$, so $a_0 s_0=0$, which implies that $a_0=0$ too, a contradiction. Hence $s_i$ are linearly independent.

Add the class of $\tilde{\Gamma}_{034}$ with $Q_0+Q_3+Q_4$. We have $\bar{s}_0$ is in the linear system of $13H-5E_0-5E_3-5E_4-7E_1-7E_2-7E_5$
$-3E_{03}-3E_{04}-3E_{34}-4E_{05}-4E_{13}-4E_{24}-2E_{12}-2E_{15}-2E_{25}$, which equals to $D-E_{12}-E_{15}-E_{25}$. The computation for the other $s_i$ is the same. As a conclusion, $h^0(X,\oO_X(D))=h^0(Y,\oO_Y(D))\geq 4$.

Now we only need to show $h^0(X,\oO_X(D))\leq 4$. Consider the restriction map $r:\Pic(X)\ra \Pic(S)$.
By Theorem \ref{pencilsQ} (1), $r(D)=\kappa^* H_S $. Hence we have the exact sequence:
\[0\ra H^0(X,\oO_X(D-S))\ra H^0(X,\oO_X(D))\ra H^0(S, \oO_S(\kappa^* H_S)).\]
Therefore $ h^0(X,\oO_X(D))\leq h^0(X,\oO_X(D-S))+h^0(S, \oO_S(\kappa^* H_S))$.
Since $\kappa$ is an automorphism of $S$, $H^0(S, \oO_S(\kappa^* H_S))=H^0(S, \oO_S(H_S))=4$. Therefore we only need to show that $D-S$ is not effective on $X$.
Let $G:=r(S)\sim -(3/2)\Lambda +(1/2)\sum_{i} E_i+\sum_{\vv{I}=3} L_I$. By (\ref{mu0}), we can compute that $\kappa^* G =G$.
We restrict $D-S$ to $S$: $r(D-mS)=\kappa^* H_S-m r(S)=\kappa^*(H_S-mG)$. 
Consider those $T_{ij}$ on $S$. We have $(H_S-mG)\cdot T_{ij}=1-2m<0$ for every $T_{ij}$ and $m\geq 1$. Now fix $m\geq 1$. Suppose now $H_S-mG$ is effective. Since each $T_{ij}$ is an irreducible $(-2)$-curve, $T_{ij}$ must lie in the fixed part of $H_S-mG$. Therefore $H_S-mG-\sum_{i\neq j} {T_{ij}}$ is effective. On the other hand,  $(H_S-mG-\sum_{i\neq j} {T_{ij}})\cdot \Lambda=6+6m-2\cdot (15)=6m-24<0$ for $m\leq 3$. Since $\Lambda$ is nef (Lemma \ref{nefbigH}), this says $H_S-mG-\sum_{i,j} {T_{ij}}$ is not effective for $m\leq 3$. Hence $H_S-mG$ is not effective for $m=1,2,3$.
As a result, $\kappa^* H_S-m G$ is not effective for $m=1,2,3$. Finally, use the exact sequences:
\[0\ra H^0(X,\oO_X(D-(m+1)S))\ra H^0(X,\oO_X(D-mS))\ra H^0(S, \oO_S(\kappa^* H_S-mG)).\]
Then we find $h^0(X,\oO_X(D-S))= h^0(X,\oO_X(D-2S))=\cdots = h^0(X,\oO_X(D-4S))=0$ where $D-4S$ is not effective because its degree is $-3<0$. This proves that $ h^0(X,\oO_X(D))=4$.\qed

\begin{defi}\label{psi}
	For fixed six points $p_0,\cdots, p_5$ in $\PP^3$, we define the rational map $\psi:\PP^3\dashrightarrow \PP^3$ by $\psi: [x_0:x_1:x_2:x_3]\mt [s_0:s_1:s_2:s_3]$.
\end{defi}

\subsection{Extra relations from quintics}
 Recall Definition \ref{DQ} that $A:= 4H-2\sum_{i=0}^5 E_i-(E_{03}+E_{04}+E_{34})-(E_{05}+E_{13}+E_{24})$. We define $6$ quintic classes in $Y$ as follows.
\[\begin{array}{ll}
	D_{05}:=&A+H-E_0-E_5-E_{03}-E_{04}-E_{05},\quad F_{15}:=A+H-E_1-E_5-E_{05}-E_{13}-E_{15},\\
	D_{13}:=&A+H-E_1-E_3-E_{03}-E_{34}-E_{13},\quad F_{25}:=A+H-E_2-E_5-E_{05}-E_{24}-E_{25},\\
	D_{24}:=&A+H-E_2-E_4-E_{04}-E_{34}-E_{24},\quad F_{12}:=A+H-E_1-E_2-E_{13}-E_{24}-E_{12}.
	\end{array}\]

	The three quintics $D_{ij}$ (and $F_{ij}$) are in the same orbit under the $\mathcal{S}_3$-action in Remark \ref{qrmk}.

	\begin{prop}\label{pencils}
		If the six points $p_0,\cdots, p_5$ are very general, then for each $D_{ij}$ and $F_{ij}$ above, $h^0(Y,\oO_Y(D_{ij}))=h^0(Y,\oO_Y(F_{ij}))=2$. Same results hold  over $X$.
\end{prop}

\begin{lemma}\label{rdual}
	In $Y$, the linear systems of $D_{ij}$ and $F_{ij}$ have the following sections:
\[\begin{array}{c|rrrr}
		D_{05} & \overline{p}_{034}\overline{f}_0 x_{E_{15}}x_{E_{25}}  & \overline{p}_{045}\overline{f}_{24}x_{E_{4}} &\overline{p}_{035}\overline{f}_{13}x_{E_{3}}\\
		D_{13} & \overline{p}_{034}\overline{f}_3 x_{E_{12}}x_{E_{15}}  & \overline{p}_{013}\overline{f}_{05}x_{E_{0}} &\overline{p}_{134}\overline{f}_{24}x_{E_{4}}\\
		D_{24} & \overline{p}_{034}\overline{f}_4 x_{E_{12}}x_{E_{25}}  & \overline{p}_{234}\overline{f}_{13}x_{E_{3}} &\overline{p}_{024}\overline{f}_{05}x_{E_{0}}\\
\end{array}\]
\label{tdual}
\[\begin{array}{c|rrrrrr}
		F_{15} & \overline{p}_{134}\overline{f}_0 x_{E_{25}}  & \overline{p}_{045}\overline{f}_3 x_{E_{12}} &  & &  \overline{p}_{015}\overline{f}_{25}x_{E_{0}} & \overline{p}_{135}\overline{f}_{12}x_{E_{3}} \\
               F_{25} & \overline{p}_{234}\overline{f}_0 x_{E_{15}} & &  \overline{p}_{035}\overline{f}_4 x_{E_{12}} & \overline{p}_{025}\overline{f}_{15}x_{E_{0}}& &\overline{p}_{245}\overline{f}_{12}x_{E_{4}} \\
	        F_{12} & &  \overline{p}_{024}\overline{f}_3 x_{E_{15}} & \overline{p}_{013}\overline{f}_4 x_{E_{25}} & \overline{p}_{123}\overline{f}_{15}x_{E_{3}}& \overline{p}_{124}\overline{f}_{25} x_{E_{4}}
\end{array}\]

\end{lemma}
\pf. We only need to prove that the sums of the divisor classes on the right equal $D_{ij}$ or $F_{ij}$. By symmetry, we only need to verify for $D_{05}$ and $F_{15}$, which follows from a direct calculation.\qed

\begin{lemma}
	In the K3 Kummer surface $S$ with $\rho(S)=17$, let $A_{ij}:=\Lambda-(1/2)(E_i+E_j+\sum_{p\neq q\in \{0,1,2,3,4,5\}-\{i,j\}} L_{ipq})$. Then $H^0(S,\oO_S(A_{ij}))=2$.
	\label{2class}
\end{lemma}
\pf. By symmetry, we only need to prove the lemma for $A_{01}$. We find $A_{01}\sim E_2+L_{012}+T_{02}+T_{12}$ is a sum of four $(-2)$-curves.
Since $A_{01}$ is effective and $A_{01}\neq \oO_S$, $h^2(S,\oO_S(A_{01}))=0$. Note that $A_{01}^2=0$. Then by Riemann-Roch: $h^0(S,\oO_S(A_{01}))\geq (1/2)A_{01}^2+2=2$.
For the other inequality we restrict $A_{01}$ to $E_2\cong \PP^1$:
\[0\ra H^0(S,\oO_S(A_{01}-E_2))\ra H^0(S,\oO_S(A_{01}))\ra H^0(E_2, \oO_{E_2}({A_{01}}_{\mid E_2})).\] 
Here the restriction map $r_2:\Pic(S)\ra \Pic(E_2)\cong \Z$ is given by  $E_2\mt -2$ and $\Lambda, E_j, L_I\mt 0$ for $j\neq 2$. Therefore by linearity, $r_2(A_{01})\sim 0$, so $H^0(E_2, \oO_{E_2}({A_{01}}_{\mid E_2}))=1$. Therefore we only need to  prove $H^0(S,\oO_S(A_{01}-E_2))\leq 1$. 
Now $A_{01}-E_2\sim T_{02}+T_{12}+L_{012}$.  We compute $( T_{02}+T_{12}+L_{012})\cdot T_{02}=-1<0$. Since $T_{02}$ is irreducible, $T_{02}$ is contained in the fixed part of  $T_{12}+L_{012}$. Therefore we only need to show $h^0(S,\oO_S(T_{12}+L_{012}))\leq 1$. Now $(T_{12}+L_{012})\cdot T_{12}=-1<0$, so $T_{12}$ is contained in the fixed part of $T_{12}+L_{012}$. Therefore we only need $h^0(S,\oO_S(L_{012}))\leq 1$, which holds since $L_{012}$ is a $(-2)$-curve. \qed

\pfof{Proposition \ref{pencils}.} By symmetry, we only need to prove the Proposition for $D_{05}$ and $F_{15}$. By Lemma \ref{rdual}, we only need to show  $h^0(X,\oO_X(D_{05}))\leq 2$ and $h^0(X,\oO_X(F_{15}))\leq 2$.
We restrict $D_{05}$ to $S$ and consider the exact sequence:
\[0\ra H^0(X,\oO_X(D_{05}-S))\ra H^0(X,\oO_X(D_{05}))\ra H^0(S, \oO_S(r({D_{05}}))).\]
Now calculation shows that $r(D_{05})=\kappa^* A_{12}$. Hence $h^0(S, \oO_S(r({D_{05}})))=h_0(S,\oO_S(A_{12}))=2$ by Lemma \ref{2class}. On the other hand, $D_{05}-S\sim H-E_0-E_5-E_{03}-E_{04}-\txt{(other } E_{ij})$ is not effective. Hence $h^0(X,\oO_X(D_{05}))\leq 2$.

Similarly, we find $r(F_{15})=\kappa^* A_{24}$, and $F_{15}-S\sim H-E_1-E_5-E_{05}-E_{13}-\txt{(other } E_{ij})$ is not effective. By a similar exact sequence, $h^0(X,\oO_X(F_{15}))\leq 2$. \qed

\begin{corollary}
	\label{extrasections}
	Let $\Span(f,g)$ be the linear span of two polynomials $f$ and $g$ over $\CC$. Then for very general choice of the six points, we have

\noindent\begin{minipage}{.5\linewidth}
\begin{align}
	\begin{split}
	 p_{034}f_0\in & \Span (p_{045}f_{24},p_{035}f_{13});\\
	 p_{034}f_3\in & \Span(p_{013}f_{05},p_{134}f_{24});\\
	 p_{034}f_4\in & \Span(p_{234}f_{13},p_{024}f_{05}).
	\end{split}
\end{align}
\end{minipage}
\begin{minipage}{.5\linewidth}
\begin{align*}
	\begin{split}
	p_{135}f_{12}, p_{015}f_{25} \in \Span(p_{134}f_0,p_{045}f_3);\\
	p_{245}f_{12}, p_{025}f_{15} \in \Span(p_{234}f_0,p_{035}f_4);\\
	p_{123}f_{15}, p_{124}f_{25} \in \Span(p_{024}f_3,p_{013}f_4).
	\end{split}
\end{align*}
\end{minipage}
\end{corollary}
Now we can define some additional degree $13$ polynomials, which give additional sections of the restriction of $\vv{D}$ to $\PP^3$. Section \ref{BIR} will show that these polynomial identities in Corollary \ref{extrasections} corresponds to the construction of the configuration $\{q_i, l_{ij}\}$ on the target.
\begin{defi}\label{sothers}
\begin{align}
	&s'_0=s_0=p_{034} f_0 f_3 f_4,\, s'_1= p_{035} f_3 f_4 f_{13},\,  s'_2=p_{024} f_0 f_3 f_{05},\, s'_3=p_{134} f_0 f_4 f_{24}.\\
	&s''_0=p_{245}f_{12} f_{3} f_{13}, \quad s''_3 =p_{135}f_{12} f_{4} f_{24}.
\end{align}
\end{defi}
\begin{corollary}
	\begin{enumerate}
		\item There exist suitable choices of the scalar multiples of $f_\alpha$ and $p_{ijk}$ such that $s_0=s_1-s'_1=s_2-s'_2=s_3-s'_3$.
\item Any of the following maps $\PP^3\dashrightarrow \PP^3$ equals $M\circ \psi$ for some $M\in \PGL(4)$:
	\begin{align*}
		&[s_0:s_1':s_2':s_3'],\\
		&[s_1':s_1:s_2:s_3],\\
		&[s_0'':s_1:s_2:s_3'']=[d s'_1+ c s_2:s_1:s_2:a s'_3+ b s_1]\\
		&\txt{ for some nonzero scalars $a,b,c,d$}.
	\end{align*}
	\end{enumerate}
	\label{salt}
\end{corollary}

\section{Birationality}\label{BIR}
In this section we prove that the rational map $\psi:\PP^3\dashrightarrow \PP^3$ is birational by constructing its inverse $\psi\ik$, which is induced by $\vv{D'}$ of a divisor class $D'$ symmetric to $D$.

\subsection{Six points on the target}
As a preparation we show that there are six special points $q_i$, $i=0,\cdots, 5$ on the target $\PP^3$ such that the quartics $Q_\alpha'$ are contracted by $\psi$ to the points $q_i$ or lines $\overline{q_i q_j}$ indexed by $\alpha$. We claim:

\begin{defiThm}	Each $Q_\alpha'$ is contracted by $\psi$ to a line or a point. In particular:
	\begin{enumerate}
		\item $Q_0'$, $Q_3'$, $Q_4'$ are contracted to the points:
			\[q_0:=[0:1:0:0],\quad q_3:=[0:0:0:1],\quad q_4:=[0:0:1:0]\]
			respectively.
		\item The quartics $Q_{05}', Q_{13}', Q_{24}'$ are contracted to three lines $l_{05}$, $l_{13}$, $l_{24}$ respectively, where $l_{05}$ passes through $q_0$, $l_{13}$ passes through $q_3$, and $l_{24}$ passes through $q_4$. 
		\item The quartics $Q_{12}', Q_{15}', Q_{25}'$ are contracted to three lines $l_{12}, l_{15}, l_{25}$ respectively.
		\item The lines $l_{12}$, $l_{15}$ and $l_{13}$ meet at a unique point $q_1$. The lines $l_{15}$, $l_{25}$ and $l_{05}$ meet at a unique point $q_5$. The lines $l_{25}$, $l_{12}$ and $l_{24}$ meet at a unique point $q_2$.
	\end{enumerate}\label{6dual}
\end{defiThm}
\pf. By symmetry, we only need to prove (1) to (3) for $Q_0', Q_{05}'$ and $Q_{12}'$.

(1). Recall that $\psi=[p_{034}f_{0}f_{3}f_{4}: p_{045}f_{3}f_{4}f_{24}:p_{234}f_{0}f_{3}f_{13}: p_{013}f_{0}f_{4}f_{05}]$. Here $f_{0}$ is the polynomial defining $Q_0'$ while $p_{045}f_{3}f_{4}f_{24}$ does not vanish on $Q_0'$. So $\psi(Q_0')=[0:1:0:0]$.
    
    (2). Recall that for suitable choice of the multiples of $f_{\alpha}$ and $p_{ijk}$, we have $s_0=s_2-s_2^{\prime}$ (Corollary \ref{salt}). 
    
Then there exists a matrix $M_{05}\in \PGL(4)$ such that $M_{05}\circ\psi=[s_0:s_1:s_2^{\prime}:s_3]$. Since $s_2^{\prime}$ and $s_3$ vanish on $Q_{05}'$ while $s_0$, $s_1$ do not, $(M_{05}\circ\psi)(Q_{05}')$ is contained in the line $\{[x:y:0:0]\mid x,y\in \CC\}$. Hence $\psi(Q_{05}')$ is contained in the line $\{[x:y:x:0]\mid x,y\in \CC\}$ which we call $l_{05}$. This line $l_{05}$  
contains the point $q_0=[0:1:0:0]$. 

(3). 
By Corollary \ref{salt}, there exists an $M_{12}\in \PGL(4)$ such that 
\[M_{12}\circ \psi=[ds_1^{\prime}+cs_2:s_1:s_2:as_3^{\prime}+bs_1]\]
for some non-zero scalars $a,b,c,d$.
    Then $M_{12}\circ\psi(Q_{12}')$ is contained in the line $\{[0:x:y:0] \mid x, y\in \mathbb{C}\}$. Hence $\psi(Q_{12}')$ is contained in the line $M_{12}^{-1}([0:x:y:0])$ which we could name as  $l_{12}$.
   
    To finish the definition, we define the lines $l_{13}$ and $l_{24}$ using the relations $s_0=s_3-s_3'=s_1-s_1'$ respectively. We define  $l_{15}$, $l_{25}$ using the polynomial identities from Corollary \ref{extrasections} similar to the one used for $l_{12}$.

    (4). We prove that $l_{12}$, $l_{24}$, $l_{25}$ meet at a unique point.
    To do this, let $\phi_D:Y\dashrightarrow \PP^3$ be the map induced by $\vv{D}$. We use Lemma \ref{E4q2} below to show that $\phi_D$ contracts $E_4$ to a point, which we defined as $q_2$. Then we prove that $q_2$ is on $l_{12}$, $l_{24}$ and $l_{25}$.
    
\underline{$q_2$ lies in $l_{12}$.} We state the following small lemma which follows from a local calculation:
    \begin{lemma}\label{smalllm}
        Let $\omega:\mathbb{P}^n\dashrightarrow\mathbb{P}^n$ be a rational map given by homogeneous polynomials $[t_0: t_1: \cdots: t_n]$. Let $W\rightarrow \mathbb{P}^n$ be the blow-up of $\mathbb{P}^n$ at a point $p$ with exceptional divisor $E_p$. Let $\tilde{\omega}: W\dashrightarrow \mathbb{P}^n$ be induced by $\omega$.  Let $x_0,\cdots, x_n$ be the coordinate functions on the target $\mathbb{P}^n$. Suppose $t_i$ vanishes at $p$ with multiplicity $m_i$. Let $m=min_i\{m_i\}$. Then for any $m_j>m$, $\tilde{\omega}(E_p)$ is contained in the coordinate plane $\{x_j=0\}$.
    \end{lemma}

    Recall the proof of (3) above that there exists an $M_{12}\in \PGL(4)$ such that $M\circ \psi=[ds_1^{\prime}+cs_2: s_1: s_2: as_3^{\prime}+bs_1]$ where $a$, $b$, $c$, $d$ are non-zero scalars. And by Corollary \ref{extrasections} we have
    \begin{align*}
        ds_1^{\prime}+cs_2=\lambda p_{245}f_{12}f_3f_{13},\\
        as_3^{\prime}+bs_1=\mu p_{135}f_{12}f_4f_{24},
    \end{align*}
    where $\lambda$ and $\mu$ are non-zero scalars. 
    Since the multiplicities of the polynomials $p_{245}$, $f_{12}$, $f_3$, $f_{13}$ at $p_4$ are 1, 2,  1, 2 respectively, the multiplicity of $ds_1^{\prime}+cs_2=\lambda p_{245}f_{12}f_3f_{13}$ at $p_4$ is 6. Similarly, we could prove that the multiplicities of $s_1$, $s_2$ and $as_3^{\prime}+bs_1$ at $p_4$ are 5, 5 and 6 respectively.
    Hence by Lemma \ref{smalllm}, when lifted to $Y$ we have
    $M_{12}\circ\phi_D(E_4)=\{[0:x_0:y_0:0]\}$ for some non-zero  $x_0, y_0 \in \mathbb{C}$. 
    Recall that $M_{12}(l_{12})$ is the line $\{[0:x:y:0] \mid x,y \in \mathbb{C}\}$, hence $M_{12}\circ \phi_D(E_4)\in M_{12}(l_{12})$. Therefore the unique point $q_2=\phi_D(E_4)$ is in $l_{12}$.  
    
    \underline{$q_2$ lies in $l_{25}$.} This is symmetric to $l_{12}$.  

\underline{$q_2$ lies in $l_{24}$.}
 Since $s_0=s_3-s_3^{\prime}$, let 
    \[M_{24}=
    \begin{bmatrix}
    1&0&0&0\\
    0&1&0&0\\
    0&0&1&0\\
    -1&0&0&1
    \end{bmatrix} \in \PGL(4),\]
    then $M_{24}\circ \psi= [s_0: s_1: s_2: s_3^{\prime}]$. Since $s_1$ and $s_3'$ vanish on $Q_{24}$ while $s_0$ and $s_2$ do not, $M_{24}\circ \phi(Q_{24}')$ is contained in the line  $\{[x:0:y:0] \mid x,y\in \mathbb{C}\}$, which we define as $l_{24}$. And multiplicities of $s_0$ and $s_2$ at $p_4$ are 5 while multiplicities of $s_1$, and $s_3^{\prime}$ at $p_4$ are 6. Therefore $M_{24}\circ \phi_D(E_4)=[x_0:0:y_0:0]$ for some $x_0, y_0\in \mathbb{C}$. Hence $M_{24}(q_2)= M_{24}\circ \phi_D(E_4)\in M_{24}(l_{24})$. So $q_2\in l_{24}$.

    Hence $l_{12}$, $l_{24}$ and $l_{25}$ intersect at the point $q_2$. The other two claims follow from symmetry.\qed
    
    \begin{lemma}
	    The map $\phi_D:Y\dashrightarrow \PP^3$ contracts $E_4$ to a point $q_2$.
	    \label{E4q2}
    \end{lemma} 

    \pf.  Recall Proposition \ref{wd} that  $h^0(Y,\oO_Y(D))=4$. To prove that $\phi_D$ contracts $E_4$ to a point we only need to prove that $h^0(Y, \oO_Y(D-E_4))=3$. Consider the exact sequence
    \[0 \rightarrow H^0(Y,  \oO_Y(D-E_4))\rightarrow H^0(Y,  \oO_Y(D))\rightarrow H^0(E_4,  \oO_Y(D_{\vert_{E_4}})).\]
    $E_4$ is isomorphic to the blow-up of $\mathbb{P}^2$ at 5 general points, say $t_0,t_1, t_2,t_3$ and $t_5$. Let $e_i$ be the exceptional divisor over $t_i$. Then 
    \[D\vert_{E_4}\sim 5l-3e_0-4e_2-3e_3.\] On $\mathbb{P}^2$ there is only one degree-5 curve whose multiplicities at  $t_0$, $t_2$ and $t_3$ are 3, 4 and 3 respectively.
    To see this, we could assume without loss of generality that $t_0=[1:0:0]$, $t_2=[0:1:0]$, $t_3=[0:0:1]$ on $\mathbb{P}^2=\proj \CC[x,y,z]$. Then we can check that up to scalar there is only one polynomial vanishing at $t_0$, $t_2$ and $t_3$ with  multiplicities $3$, $4$ and $3$ respectively, which is $x^2yz^2$.
    Hence $h^0(D\vert_{E_4})=1$. So by the exact sequence, 
    $h^0(D-E_4)\geq 3$. The section $s_0 x_{E_{12}} x_{E_{15}}x_{E_{25}}$ of $D$ is not in the image of $H^0(D-E_4)$. Hence $h^0(D-E_4)\leq 3$. Hence $h^0(D-E_4)=3$.\qed
    
\begin{figure}[ht]
	\centering
	\includegraphics[height=5cm,clip]{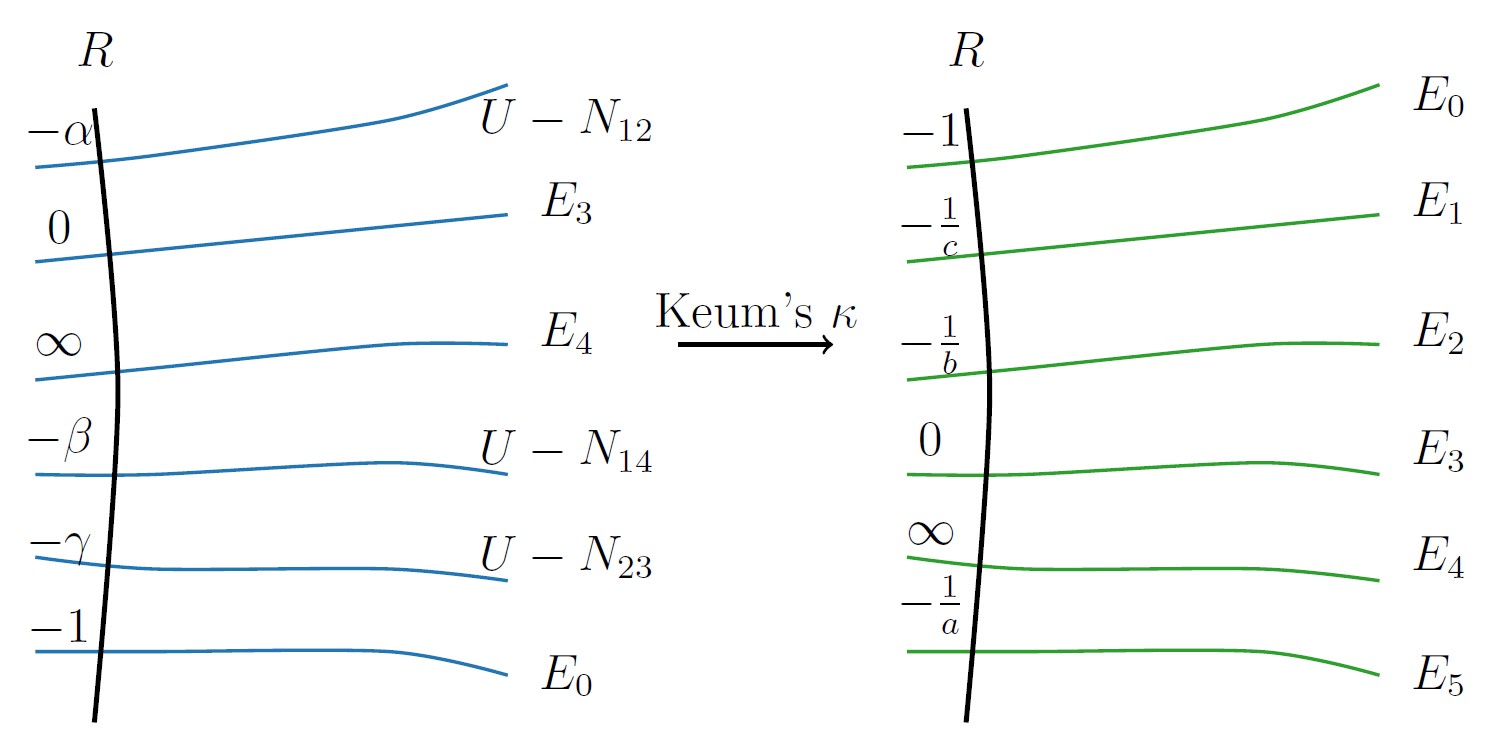}
	\caption{The six points on the rational normal curve $R$ as intersections with the $(-2)$ curves}
	\label{fig:sixpoints}
\end{figure}

    \begin{theorem}If the six points $p_i$ are in very general position, then the six points $q_i$, $i=0,\cdots, 5$ are distinct, and there exists an $M\in \PGL(4)$ such that $M$ sends $(q_0,\cdots, q_5)$ to $(p_0, \cdots, p_5)$.  That is, the six points $q_i$ are projectively equivalent to $p_i$. \label{projequiv}
\end{theorem}
\pf.  Without loss of generality, we can assume the six general points on $\mathbb{P}^3$ to be
    \begin{align*}
	    &p_0=[1:0:0:0],  \quad  p_1=[0:0:0:1],  \quad p_2=[0:0:1:0],\\
	    & p_3=[1:1:1:1], \quad  p_4=\left[1:\frac{1}{a}:\frac{1}{b}:\frac{1}{c}\right], \quad p_5=[0:1:0:0],
    \end{align*}
    for $a,b,c$ nonzero.
    Then the rational normal curve $R_0$ in $\PP^3$ through these six points can be given by $R_0: \mathbb{P}^1\rightarrow \mathbb{P}^3$, sending $[u:v]$ to 
    \begin{align}
    &\left[\frac{1}{u+v}: \frac{1}{au+v}:\frac{1}{bu+v}:\frac{1}{cu+v}\right]\\
    =&[(au+v)(bu+v)(cu+v):\cdots:(u+v)(au+v)(bu+v)].\label{r0}
    \end{align}
Here Keum's automorphism $\kappa$ maps the following $(-2)$-curves to $E_i$ and $R$ (see Figure \ref{fig:sixpoints}):
\begin{align*}
	\{ U-N_{12},U-N_{14},U-N_{23}\}&\mapsto \{E_0,E_{3},E_4\},\\
	\{E_3,E_4,E_0\}&\mapsto \{E_1,E_2,E_5\},\\
        R&\mapsto R.
\end{align*}
 
Each $E_i$ meets $R$ at a unique point $p_i'$ over $p_i$. Therefore after applying $\kappa\ik$ we find that $U-N_{12}$, $U-N_{14}$ and $U-N_{23}$ each meets $R$ at a unique point.
The $[u:v]$ coordinates of $p_i$ in $R_0\subseteq \mathbb{P}^3$ are equal to the $[u:v]$ coordinates of $p_i^{\prime}$ in $R$ because $R_0\cong R$.
Hence there exists $\alpha, \beta, \gamma \in \CC\cup \{\infty\}$ such that the $[u:v]$ coordinates for the six points where $U-N_{12}$, $E_3$, $E_4$, $U-N_{14}$, $U-N_{23}$, $E_0$ meet $R$ are given by
\[u/v=\{-\alpha,0,\infty,-\beta,-\gamma,-1\},\]
Since $\kappa$ maps $R$ isomorphically to $R$ itself, $\kappa_{\vert_{R\cong \mathbb{P}^1}}$ is some $M\in \PGL(2)$ which sends the six points with $[u:v]$ coordinates $(-\alpha,0,\infty,-\beta,-\gamma,1)$ to $(p_0',\cdots,p_5')$, that is, $(p_0,\cdots,p_5)$ when considered in $R_0$. 
The coordinates of $p_i$ give the $[u:v]$ coordinates of $p_0,\cdots,p_5$ on $R$:
\[u/v=\{-1,-1/c,-1/b,0,\infty,-1/a\}.\]
The cross ratio of $p_0,p_1,p_2,p_5$ equals to the cross ratio of their preimage under $\kappa$. Therefore we obtain $\alpha=(c-1)(b-a)/((b-1)(c-a))$. This implies that for $a,b,c$ general, $\alpha\not\in\{-1,-1/c,-1/b,0,\infty,-1/a\}$, which implies that the point where $U-N_{12}$ meets $R$ is not one of $p_i'$. By symmetry, this also holds for $U-N_{14}$ and $U-N_{23}$.

    We next compute $\psi(R_0)$ by restricting the polynomials $f_{\alpha}$ and $p_{ijk}$ to $R_0$.
  Define $f_{\alpha, R}$ and $p_{ijk,R}$ to be the polynomial obtained by plugging in the  $[x_0:x_1:x_2:x_3]$ in $f_{\alpha}(x_0,x_1,x_2,x_3)$ or $p_{ijk}(x_0,x_1,x_2,x_3)$ by (\ref{r0}).
  We claim that by counting multiplicities at $p_0,\cdots,p_5$, up to scalars:
  \begin{align*}
      f_{12,R}&=f_{15,R}=f_{25,R}=f_{05,R}=f_{13,R}=f_{24,R}\\
      &=u^2v^2(u+v)^2(au+v)^2(bu+v)^2(cu+v)^2,\\
      f_{0,R}&= uv(u+v)^2(au+v)^3(bu+v)^2(cu+v)^2(u+\alpha v),\\
      f_{3,R}&= u^2v(u+v)(au+v)^2(bu+v)^2(cu+v)^3(u+\beta v),\\
      f_{4,R}&= uv^2(u+v)(au+v)^2(bu+v)^3(cu+v)^2(u+\gamma v),\\
 	p_{034,R}&=uv(u+v),   \\ 
	p_{045,R}&=v(u+v)(au+v), \quad  	p_{234,R}=uv(bu+v), \quad p_{013,R}=u(cu+v)(u+v).
  \end{align*}

  Here the only nontrivial parts are the terms containing $\alpha, \beta$ and $\gamma$. We prove the equality for $f_0$ above, and by symmetry it holds for $f_3$ and $f_4$.  Recall Theorem \ref{pencilsQ} (1) that $Q_0$ restricts to the $(-2)$-curve $\kappa^*E_0=U-N_{12}$ on $S$. Therefore in $X$ we have $Q_0\cap R=Q_0\cap S\cap R=(U-N_{12})\cap R$ is a unique point $p_0''$ whose $[u:v]$ coordinate is $[-\alpha:1]$. Thus $(u+\alpha v)\mid f_{0,R}$. Since $\alpha\not\in\{-1,-1/c,-1/b,0,\infty,-1/a\}$, we know $p_0''\neq p_i'$ for every $i$ when $a,b,c$ are general. Hence $(u+\alpha v)$ is not one of the factors in $ uv(u+v)^2(au+v)^3(bu+v)^2(cu+v)^2$.  By counting multiplicities, $f_{0,R}$ equals the product on the right-hand side.

Now we know the image $R_0$ under $\psi$ is given by
\begin{align*}
	R'([u:v]):=\psi(R_0([u:v]))&=\left[p_{034,R} :p_{045,R}\frac{f_{24,R}}{f_{0,R}}: p_{234,R} \frac{f_{13,R}}{f_{4,R}}: p_{013,R}\frac{f_{05,R}}{f_{3,R}}\right]\\
&=\left[1:\frac{ v}{u+\alpha v}:\frac{u}{u+\gamma v}:\frac{u+v}{u+\beta v}\right]\in \mathbb{P}^3.
\end{align*}

As a result, $R'$ is a degree-3 rational curve.
Next we prove that the six points $q_i$ are all on $R^{\prime}$ and we find the $[u:v]$ parameter for them in $\mathbb{P}^1$.
Note that $q_0=[0:1:0:0]=R^{\prime}([-\alpha:1])$, $q_3=[0:0:0:1]=R^{\prime}([-\beta:1])$, and $q_4=[0:0:1:0]=R^{\prime}([-\gamma:1])$.
By Lemma \ref{E4q2}, $\phi_D:Y\dashrightarrow \PP^3$ contracts $E_3$, $E_4$, $E_0$  to $q_1$, $q_2$, $q_5$ respectively.  Note that $E_3\cap R=R([0:1])$, $E_4\cap R=R([1:0])$, and $E_0\cap R=R([1:-1])$. We claim that
\begin{align*}
q_1&=\phi_D (E_3)=R'([0: 1]),\quad q_2=\phi_D (E_4)=R'([1:0]),\quad q_5=\phi_D (E_0)=R'([1:-1]).
\end{align*}
Indeed, by symmetry we only need to verify that $\phi_D$ is defined at the point $p_3'$ where $E_3$ intersects $R$. To see this, we only need to show $\bar{s}_1(p_3')\neq 0$. Here $s_1=p_{045}f_3 f_4 f_{24}$. Clearly $p_{045}$ does not vanish at $p_3'$. Then $Q_3$ and $Q_4$  restrict to $U-N_{14}$ and $U-N_{23}$, each of which intersects $R$ at a different point from $p_3'$ by the proof above (Figure (\ref{fig:sixpoints})). Finally, $Q_{24}$ restricts to $\kappa\ik(T_{24})$, where $\kappa\ik(T_{24})\cap R=\kappa\ik(T_{24}\cap R)$ is empty. Together we find $\bar{s}_1(p_3')\neq 0$.

Now we know that the six points $q_0,q_1,\cdots,q_5$ are on $R'$, corresponding to the six numbers 
\[u/v=\{-\alpha, 0, \infty,-\beta,-\gamma,-1\}.\]
Therefore  $(q_0,\cdots,q_5)$ are distinct, and $R'$ is the unique rational normal curve through $q_0,\cdots, q_5$.
As a result, the matrix $M=\kappa_{\mid R}\in \PGL(2)$ sends $(q_0, \cdots, q_5)$ to $(p_0,\cdots,p_5)$. This implies that $(p_0,\ldots,p_5)$ is projectively equivalent to $(q_0,\ldots,q_5)$.\qed

\subsection{The dual construction}
As a corollary of Theorem \ref{projequiv}, $l_{ij}$ in Definition-Theorem \ref{6dual} equals the line $\overline{q_i q_j}$. Hence the six distinct points $q_i$ and the $9$ lines $l_{ij}$ on the target form the same configuration as $p_i$ and $l_{ij}$ for $ (ij)\in \mathcal{I}$ in the source. Hence blowing up the $6$ points and $9$ lines in the source and target induces a rational map $\phi:Y\dashrightarrow Y$.

We now define divisor classes $D'$ and $9$ dual quartics classes $P_\beta$ on the target, where $\beta\in \mathcal{B}:=\{1,2,5,03,04,34,05,13,24\}$, by switching the index $1$ with $3$, $2$ with $4$ and $0$ with $5$ in the classes of $D$ and $Q_\alpha$. That is:
\begin{align*}
D':=\quad &  13H-5(E_1+E_2+E_5)-7(E_0+E_3+E_4)\\
	&-(E_{03}+E_{04}+E_{34})-4(E_{05}+E_{13}+E_{24})-3(E_{12}+E_{15}+E_{25}).
	\end{align*}
\[
	\begin{array}[ht]{lll}
	P_1 :=& 4H-2E_1-E_2-E_5-2E_0-3E_3-2E_4\\
	&\hspace{5em}-E_{12}-E_{15}-E_{05}-2E_{13}-E_{24}-E_{03}-E_{04},\\
	P_2 :=& 4H-E_1-2E_2-E_5-2E_0-2E_3-3E_4\\
	&\hspace{5em}-E_{12}-E_{25}-E_{05}-E_{13}-2E_{24}-E_{03}-E_{34},\\
	P_5 :=& 4H-E_1-E_2-2E_5-3E_0-2E_3-2E_4\\
	&\hspace{5em}-E_{15}-E_{25}-2E_{05}-E_{13}-E_{24}-E_{04}-E_{34},\\
	P_{05} :=& B-E_{12},\qquad P_{13} := B-E_{25}, \qquad P_{24} := B-E_{15},\\
	P_{34} :=& B-E_{05},\qquad P_{03} := B-E_{24}, \qquad P_{04} := B-E_{13},\\
	\end{array}
\]
where $B:=4H-2\sum_{i=0}^5 E_i-(E_{12}+E_{15}+E_{25})-(E_{05}+E_{13}+E_{24})$.

Similar to $Q_\alpha$, we can define $P_\beta'$ as the image of $P_\beta$ in $\PP^3$. Let $g_\beta$ be the polynomial defining $P_\beta'$.  Let $q_{ijk}$ be the polynomial defining the plane in $\PP^3$ through the points $q_i,q_j$ and $q_k$. Now we define
:
\begin{align}
	\begin{split}
	\psi':\quad  & \PP^3\dashrightarrow \PP^3\\
	&[y_0:y_1:y_2:y_3]\mt [t_0:t_1:t_2:t_3]
	\end{split}
\end{align}
with
\[
[t_0:t_1:t_2:t_3]:=
[q_{125} g_1 g_2 g_5: q_{025} g_1 g_2 g_{24}:  q_{124} g_1 g_5 g_{13}: q_{135} g_2 g_5 g_{05}].\]
The symmetry between the six points $p_i$ and $q_i$ implies that those $P_\beta$ satisfy the dual version of Theorem \ref{pencilsQ}. The map $\psi'$ is induced by $\vv{D'}$. 

We introduce some notations. For each $\alpha \in \mathcal{A}$ and  $\beta\in \mathcal{B}$, let $m^{\alpha}_\beta$ be the multiplicity of the exceptional divisor $E_\alpha$ in the class $P_\beta$. For instance, $m^{0}_{1}=2$, and $m^{05}_{03}=1$, and $m^{03}_{2}=0$. We claim:
\begin{prop}\label{fusion}
	For each  $\beta\in \mathcal{B}$, the composition $g_\beta(\psi)=g_\beta(s_0,s_1,s_2,s_3)$ is a degree-$52$ polynomial, which up to a nonzero scalar is a product of $f_\alpha$, $\alpha\in \mathcal{A}$, in the following way:
	\begin{equation}
		g_\beta (\psi)=\prod_{\alpha\in \mathcal{A}} f_\alpha^{m^{\alpha}_{\beta}}.
	\end{equation}
\end{prop}
\pf. We first show $g_\beta (\psi)$ is not a zero polynomial. This follows from the first half of the proof of Proposition \ref{E12image}, where we show that the $\phi_D(E_{12})=l_{34}$. We note that the whole proof of  Proposition \ref{E12image} is local at $E_{12}$ and does not require that $\psi$ is birational. Therefore the image of $\psi$ contains $l_{34}$, and by symmetry, the lines $l_{03}$ and $l_{04}$. Here no $P_\beta$ passes through all the three lines $l_{34},l_{03}$ and $l_{04}$, hence the image of $\psi$ is not contained in $P_\beta$, so $g_\beta (\psi)\not\equiv 0$.  

Let  $Z_\alpha$ be the point or line indexed by $\alpha$ in the target. Now if $\mathbf{x}\in Q_\alpha$, then $\psi(\mathbf{x})$ is contained in $Z_\alpha$ by Definition-Theorem \ref{6dual}. If  $m^{\alpha}_{\beta}>0$, then $P_\beta$ passes through $Z_\alpha$, so $g_{\beta}(\psi(\mathbf{x}))=0$. Since $f_\alpha$ is irreducible, we have $f_\alpha\mid g_\beta (\psi)$. We claim the multiplicity of $f_\alpha$ in the composition $g_\beta(\psi)$ is at least $m^{\alpha}_\beta$. 
Indeed, in general let $f:U_1\ra U_2$ be a morphism. Suppose $Z_1$ is a closed subvariety in $U_1$ and $f(Z_1)=Z_2$  a closed subvariety of $U_2$. Suppose $x_1$ is the general point of $Z_1$ and let $x_2:=f(x_1)$. Let $\sigma:\oO_{U_2,x_2}\ra \oO_{U_1,x_1}$ be the induced map on local rings. Then $\sigma(\mathfrak{m}_2)\subseteq\mathfrak{m}_1$, where $\mathfrak{m}_i$ is the maximal ideal of $\oO_{U_i,x_i}$. Thus $\sigma(\mathfrak{m}_2^d)\subseteq\mathfrak{m}_1^d$ for any $d\geq 1$. Now if $h$ is a regular function on $U_2$ with multiplicity $d$ at $Z_2$, then $h\in \mathfrak{m}_2^d$, so  $f\circ h\in \mathfrak{m}_1^d$, that is, $f\circ h$ has  multiplicity at least $d$ at $Z_1$.

Run this for all $\alpha$, we find $G:= \prod_{\alpha\in \mathcal{A}} f_\alpha^{m^{\alpha}_{\beta}}$ divides $g_\beta$. Now for each $\beta$, adding up the multiplicities in $Q_\beta$ from $\mathcal{A}$ gives exactly $13$. That is:
\[\sum_{\alpha\in \mathcal{A}} m^{\alpha}_{\beta}=13.\]
Hence both $G$ and $g_\beta(\psi)$ have degree $52$. Therefore they differ by a nonzero constant. \qed

\begin{prop}
	\label{planefusion}For the $4$ planes in the definition of $\psi'$, we have up to nonzero scalars:
	\begin{align*}
		&q_{125}(\psi)=p_{125} f_{12} f_{15} f_{25}, \quad q_{025}(\psi)= p_{124} f_0 f_{05} f_{25},\\
		&q_{124}(\psi)= p_{135} f_4 f_{24} f_{12}, \quad q_{135}(\psi)= p_{025} f_3 f_{13} f_{15}.
	\end{align*}
\end{prop}
\pf. First look at $q_{125}$. Recall $\psi=[s_0:s_1:s_2:s_3]$ and Proposition \ref{wd}. Since $q_{125}$ is a plane, $q_{125}(\psi)$ is a linear combination of $s_0$, $s_1$, $s_2$, $s_3$. 
Same as above, $q_{125}(\psi)$ is not zero polynomial since the image of $\psi$ contains the line $l_{34}$, while $q_{125}$ does not vanish on $l_{34}$.
Using the same argument in the proof of Proposition \ref{fusion}, we can show that $f_{12}$, $f_{15}$, $f_{25}$ are all irreducible factors of $q_{125}(\psi)$. Therefore $q_{125}(\psi)=h f_{12}f_{15}f_{25}$ where $h$ is a linear polynomial. Since the class of each $s_i$ has the term $-E_{12}$, each of $s_i$ vanishes at the line $\overline{p_1 p_2}$, so $q_{125}(\psi)$ vanishes at $\overline{p_1 p_2}$.  Since none of $f_{12}$, $f_{15}$ and  $f_{25}$ vanishes on the line $\overline{p_1 p_2}$, $h$ must vanish on $\overline{p_1 p_2}$. Similarly, $h$ must vanish on the lines $\overline{p_1 p_5}$. Hence up to a scalar $h=p_{125}$. 

The remaining equalities follow similarly, noticing that for $q_{025}(\psi)$ we only need to verify that $q_{025}(\psi)$ vanishes at $\overline{p_1 p_2}$ and $\overline{p_2 p_4}$ with multiplicity $1$ and $4$ (directly read from the divisor classes of $s_i$), but $f_0 f_{05} f_{25}$ does not vanish at $\overline{p_1 p_2}$ and vanishes with multiplicity exactly $3$ at $\overline{p_2 p_4}$. \qed

\begin{theorem}
	There exists a matrix $M\in \PGL(4)$ whose rows are given by the coefficients of $x_i$ in $ p_{125}$, $ p_{124}$, $ p_{135}$, and $ p_{025}$. Then up to scalars
	$ \psi'\circ \psi=M\in \PGL(4)$.
	Furthermore, $\psi$ and $\psi'$ are birational maps.
	\label{m}
\end{theorem}
\pf. By Proposition \ref{fusion} and \ref{planefusion} above, we can compute $t_i(\psi)$ using $p_{ijk}(\psi)$ and $g_\beta(\psi)$. It is easy to verify that up to nonzero scalars
\[(t_0,t_1,t_2,t_3)(\phi)=( p_{125}F,   p_{124} F,  p_{135} F,  p_{025} F),\]
for 
\[F= (f_0 f_3 f_4 )^7 (f_{05} f_{13}f_{24})^4 ( f_{12} f_{15} f_{25})^3,\]
 a polynomial of degree $(7+4+3)\cdot 4\cdot 3=168$. 
 Hence canceling $F$ gives $\psi'\circ\psi=[ p_{125}:  p_{124}:  p_{135}:  p_{025}]$, which equals to $M$. Then we only need to show $M$ is nonsingular. Indeed, if we place $p_i$ at the following position:
 \[(p_0,\cdots,p_5)=([1:a : 1: 0], [0: 0: 0: 1], [0: 0: 1: 0], [1: 1:  0: b ], [1: 0: c: 1], [0: 1: 0: 0]),\]
 then $p_0,\cdots,p_5$ are in linearly general position for $a,b,c$ general.  Then  up to nonzero scalars $p_{125}=x_1, p_{124}=x_2, p_{135}=x_3$ and $p_{025}=x_4$, hence are exactly the coordinate hyperplanes. As a result, $M$ can be chosen as the identity matrix.
 Finally, all the results above hold  by symmetry if we switch $\psi'$ with $\psi$. Hence $\psi\circ \psi'=M'$ for another $M'\in \PGL(4)$. This proves that $\psi$ and $\psi'$ are birational.\qed

\subsection{The exceptional set}

We prove that the exceptional set of $\psi$ consists of exactly the nine quartics $Q_\alpha'$:
\begin{prop}
	Consider the six points $p_i$ in very general position. Let $J$ be the Jacobian matrix of $\psi: [x_0:x_1:x_2:x_3]\mt [s_0:s_1:s_2:s_3]$. Then up to a nonzero scalar, \[\det J= ( f_0^2 f_3^2 f_4^2 ) (f_{05} f_{13}f_{24}) ( f_{12} f_{15} f_{25}).\] In particular, the only hypersurfaces contracted by $\psi$ are the $9$ quartics $Q_\alpha'$, $\alpha\in \mathcal{A}$.
	\label{exc}
\end{prop}

\begin{lemma}
	Suppose $g$ and $h_i$, $i=1,\cdots, n$ are polynomials of $n$ variables $x_1,\cdots, x_n$ and $g\neq 0$. Write $J(h_1,\cdots, h_n)$ as the Jacobian of $h_1,\cdots, h_n$ with respect to $x_1,\cdots, x_n$. If $g\mid h_1,\cdots, h_m$ for some $m$ with $2\leq m\leq n$, then $g^{m-1}\mid \det J(h_1,h_2,\cdots,h_n)$.
	\label{detNLemma}
\end{lemma}
\pf.
By assumption, we can write $h_i=gf_i$ for some polynomial $f_i$, for each $1\leq i\leq m$. Define $\vr{h}_{x}:=(\partial{h}/\partial{x_1},\cdots,\partial{h}/\partial{x_n})^T$. Then
\begin{align*}
	\det J(h_1,h_2,\cdots,h_n)	&=\det\,  [\vr{gf_1}_{x},\cdots,\vr{gf_{m}}_{x},\vr{h_{m+1}}_{x},\cdots,\vr{h_n}_{x}]\\
	&=\det\,  [g\vr{f_1}_{x}+f_1\vr{g}_x,\cdots, g\vr{f_m}_{x}+f_m\vr{g}_x, \vr{h_{m+1}}_{x},\cdots,\vr{h_n}_{x}].
\end{align*}
Now expand the columns $g\vr{f_i}_{x}+f_i\vr{g}_x$ in the last expression so that $\det J(h_1,h_2,\cdots,h_n)$ equals the sum of $2^m$ determinants. If any one of these determinants contains two columns $f_i\vr{g}_x$ and $f_j\vr{g}_x$, then it equals zero. Therefore $\det J(h_1,h_2,\cdots,h_n)$ equals
\begin{align*}
	& \sum_{i=1}^m \det\, [ g\vr{f_1}_x,\cdots,  g\vr{f_{i-1}}_x, f_{i}\vr{g}_{x}, g\vr{f_{i+1}}_x,\cdots, g\vr{f_m}_x, \vr{h_{m+1}}_{x},\cdots,\vr{h_n}_{x}]\\
	+& \sum_{i=1}^m \det\, [g\vr{f_1}_x,\cdots,g\vr{f_{m}}_x, \vr{h_{m+1}}_{x},\cdots,\vr{h_n}_{x}].
\end{align*}
Now $g^{m-1}$ divides both terms above, hence $g^{m-1}\mid \det J(h_1,h_2,\cdots,h_n)$.\qed

\pfof{Proposition \ref{exc}}. By Theorem \ref{m}, $\psi$ is birational, hence the Jacobian determinant $\Delta:=\det J$ of $\psi$ is nonzero.  Up to a nonzero scalar, $\Delta$ is invariant under change of coordinates on the target. By symmetry and Lemma \ref{detNLemma}, we only need to show $f_0^2, f_{13}$ and $f_{12}$ divide $\Delta$. Then $\Delta':=( f_0^2 f_3^2 f_4^2 ) (f_{05} f_{13}f_{24}) ( f_{12} f_{15} f_{25})$ divides $\Delta$. Now by definition, $\deg\Delta\leq (13-1)\cdot (3+1)=48$, which equals the degree of the product $\Delta'$. Hence $\Delta=\lambda \Delta'$ for some nonzero scalar $\lambda$.

So we prove $f_0^2, f_{13}$ and $f_{12}$ divide $\Delta$. First, $\psi$ is defined as $[x_0:x_1:x_2:x_3]\mt [s_0:s_1:s_2:s_3]$, where $f_0\mid s_0,s_2$ and $s_3$, by Lemma \ref{detNLemma}, $f_0^2\mid \Delta$.
Next, by Corollary \ref{salt}, we can replace $\psi$ by the map $[x_0:x_1:x_2:x_3]\mt [s_1':s_1:s_2:s_3]$. We have $f_{13}\mid s_2$ and  $f_{13}\mid s'_1$. Hence by Lemma \ref{detNLemma}, $f_{13}\mid \Delta$.
Finally by Corollary \ref{salt},  we can replace $\psi$ by the map $[x_0:x_1:x_2:x_3]\mt [s_0'':s_1:s_2:s_3'']=[d s'_1+ c s_2:s_1:s_2:a s'_3+ b s_1]$ for some nonzero scalars $a,b,c,d$. Now $f_{12}\mid s''_0$ and $ s''_3$. Hence by Lemma \ref{detNLemma}, $f_{12}\mid \Delta$.
\qed

\section{Images of the Quartics}\label{imageQ}

In this section we consider $\psi:\PP^3\dashrightarrow \PP^3$ in Definition \ref{psi}. Blowing up the six points $q_i$, $i=0,\cdots, 5$ and the $9$ lines through them indexed by $\{12,15,25,03,04,34,05,13,24\}$ in the target induces a birational map $\phi:Y \dashrightarrow Y$. We show that $\phi$ does not contract any of the nine $Q_\alpha$. 
\subsection{Lemmas on Jacobian determinants} We prove some results on Jacobian determinants which we use in the next paragraph. 
\begin{lemma}\label{fakeJac}
	Suppose $h_i$ are homogeneous in $x_0,x_1,x_2,x_3$ of degree $d\geq 1$. Write the partial derivatives $\partial h_i/\partial x_j$ as $(h_i)_{x_j}$. Then
\[\Phi:=\det
	\begin{bmatrix}
		h_0 & h_1 & h_2 & h_3 \\
		(h_0)_{x_1} & (h_1)_{x_1} & (h_2)_{x_1} & (h_3)_{x_1} \\
		(h_0)_{x_2} & (h_1)_{x_2} & (h_2)_{x_2} & (h_3)_{x_2} \\	
		(h_0)_{x_3} & (h_1)_{x_3} & (h_2)_{x_3} & (h_3)_{x_3} \\
\end{bmatrix}= \left(\frac{x_0}{d}\right)\det J(h_0,h_1,h_2,h_3)_{x_0,x_1,x_2,x_3}.
\]	\end{lemma}
\pf. For each $h_i$ we have
\begin{align}
	\label{partialsum}
	h_i(x_1,\cdots,x_n)=\frac{1}{d}\sum_{j=1}^n x_j (h_i)_{x_j}(x_1,\cdots,x_n).
\end{align}
Then we can expand $\Phi$ into the weighted sum of four determinants, among which only the one with $(h_i)_{x_0}(x_0,x_1,x_2,x_3)$ is nonzero. Hence the Lemma holds. \qed

\begin{lemma} 
	\label{imagesJac}
	Suppose $h_i$ are homogeneous in $x_0,x_1,x_2,x_3$ of degree $d\geq 1$. Then
	\begin{align*}
		\det J(h_0/h_1, h_2/h_0, h_3/h_0)_{x_1,x_2,x_3}=&\frac{-x_0}{(h_0 h_1)^2 d}\det J(h_0,h_1,h_2,h_3)_{x_0,x_1,x_2,x_3}.\\
		\det J(h_1/h_0, h_2/h_0, h_3/h_2)_{x_1,x_2,x_3}=&\frac{x_0}{(h_0^3 h_2) d}\det J(h_0,h_1,h_2,h_3)_{x_0,x_1,x_2,x_3}.
	\end{align*}
\end{lemma}
\pf. We prove the first equation and the second follows from a similar argument. Write $\vr{h}_{x}:=((h)_{x_1},(h)_{x_2},(h)_{x_3})^T$. Using the quotient rule we find
\begin{align*}
	&\det J(h_0/h_1, h_2/h_0, h_3/h_0)_{x_1,x_2,x_3}\\
	=&(h_1^{-2} h_0^{-4})\det [ h_1 \vr{h_0}_x-h_0 \vr{h_1}_x,h_0 \vr{h_2}_x-h_2 \vr{h_0}_x,h_0 \vr{h_3}_x-h_3 \vr{h_0}_x]\\
	=&(h_1^{-2} h_0^{-4}) \left(h_0^2 h_1\det J(h_0,h_2,h_3) -h_0^3 \det J(h_1,h_2,h_3)\right. \\
	 & \left.\quad +h_0^2 h_3\det J(h_1,h_2,h_0) +h_0^2 h_2\det J(h_1,h_0,h_3)\right)\\
	 = & -(h_0 h_1)^{-2} \Phi.
\end{align*}
Hence the result follows from Lemma \ref{fakeJac}.\qed
\subsection{Images of $Q_\alpha$}
We abuse notations here and denote by $Q_\alpha$  the $9$ singular quartics in $\PP^3$.
\begin{prop}\label{QtoE}
	For the six points $p_1,\cdots,p_5$ in very general position,
	let $Y_{\alpha}$ be the blow-up of $\PP^3$ at the point $q_\alpha$ for $\alpha=0,3,4$ or the line $l_\alpha$ for $\alpha=05,13,24,12,15,25$. Then the lift $\psi_\alpha:\PP^3\dashrightarrow Y_\alpha$ of $\psi$ does not contract the quartic $Q_{\alpha}$.
\end{prop}
We show that Proposition \ref{QtoE} implies that $\phi$ contracts none of the quartics $Q_\alpha$. Indeed, since the blow-up $\pi:Y\ra \PP^3$ factors as $Y\dashrightarrow Y_\alpha\ra \PP^3$, the birational map $\phi:Y\dashrightarrow Y$ does not contract $Q_\alpha$. 

\pf. By symmetry, we need only prove for $\alpha=0,05$ and $12$.

\underline{Case I: $\alpha=0$.}
Recall Definition-Theorem  \ref{6dual}(1) that $\psi(Q_0)=\{q_0\}$ with $q_0=[0:1:0:0]$. Let the homogeneous coordinates on the target copy of $\PP^3$ be $[y_0:y_1:y_2:y_3]$. Then we take the local chart $U$ at $q_0$:
\[U=\left\{((y_0,y_2,y_3),[a:b:c])\Bigm|by_0 =a y_2, cy_0 =a y_3, by_3=cy_2\right\}\subset \mathbb{A}^3\times \PP^2.\]
Then take $V$ open in $U$ defined by $V:=\{a\neq 0\} =\{a=1\}$. Then $V\cong \spec\CC[y_0,b,c]\cong \mathbb{A}^3$. Here $\phi_0:\PP^3\dashrightarrow U$ is given by:
\begin{align}
	x=[x_0:x_1:x_2:x_3]\mt \left(\left(\frac{s_0}{s_1},\frac{s_2}{s_1},\frac{s_3}{s_1}\right),\left[\frac{s_0}{f_0}:\frac{s_2}{f_0}:\frac{s_3}{f_0}\right]\right).
	\label{liftingformula0}
\end{align}
On the source let $W:=\{x_0=1\}\cong \mathbb{A}^3$. Then $\psi_0:W\dashrightarrow V$ is given by a rational map $\xi_0:\mathbb{A}^3\dashrightarrow \mathbb{A}^3$, where 
\[\xi_0(x)= (y_0, b,c)=\left(\frac{s_0}{s_1},\frac{s_2}{s_0},\frac{s_3}{s_0}\right).\]
By Lemma \ref{imagesJac} and Proposition \ref{exc}, up to a nonzero scalar the Jacobian determinant $\det J(\xi_0)$ of $\xi_0$ equals
\begin{align}
	-\frac{x_0}{13(s_0 s_1)^2}\det J(s_0,s_1,s_2,s_3)_{x_0,x_1,x_2,x_3}= -\frac{x_0  f_{05}f_{13}f_{12}f_{15}f_{25} }{13 (p_{034} p_{045} f_3 f_4)^2 f_{24} }.
\end{align}
Therefore, $\det J(\xi_0)$ does not vanish at a general point in $Q_0\cap W$. Hence $Q_0$ is not contracted by $\psi_0$.

\underline{Case II: $\alpha=05$.}
By Corollary \ref{salt}, there exists some $M_{05}\in \PGL(4)$ such that $M_{05}\circ \psi$ is given by $[s_0:s_1:s_2':s_3]=[ p_{034} f_0 f_3 f_4:  p_{045} f_3 f_4 f_{24}:   p_{024} f_0 f_3 f_{05}: p_{013} f_0 f_4 f_{05}]$.
Therefore we prove the same statement for $M_{05}\circ \psi$, where $Q_{05}$ is contracted to the line $l_{05}=\{[*:*:0:0]\}$. We take the open $U$ containing $l_{05}\backslash\{[0:1:0:0]\}$ given by
\[U:=\left\{(y_1,y_2,y_3),[a:b]\Bigm|by_2 =ay_3\right\} \subset \mathbb{A}^3\times \PP^1.\]
Take the open $V:=\{a=1\}\subset U$. Then $V$ is affine: $V\cong \spec \CC[y_1,y_2,b]\cong \mathbb{A}^3$.
Now $\psi_{05}:\PP^3\dashrightarrow U$ is defined by
\begin{align}
	x=[x_0:x_1:x_2:x_3]\mt \left(\left(\frac{s_1}{s_0},\frac{s_2'}{s_0},\frac{s_3}{s_0}\right),\left[\frac{s_2'}{f_0 f_{05}}:\frac{s_3}{f_0 f_{05}}\right]\right).
	\label{liftingformula05}
\end{align}
Then locally on $W$, $\psi_{05}$ is given by the rational map $\xi_{05}:\mathbb{A}^3\dashrightarrow \mathbb{A}^3$, where
\[\xi_{05}(x)= (y_1, y_2,b)=\left(\frac{s_1}{s_0},\frac{s_2'}{s_0},\frac{s_3}{s_2'}\right).\]
Therefore by Lemma \ref{imagesJac} and Proposition \ref{exc}, up to a nonzero scalar, 
\begin{align*}
\det J(\xi_{05})=&\frac{x_0}{s_0^3 s_2'}\det J(s_0,s_1,s_2',s_3)_{x_0,x_1,x_2,x_3}\\
=&\frac{x_0}{s_0^3 s_2'}\det J(s_0,s_1,s_2,s_3)_{x_0,x_1,x_2,x_3} \\
=&\frac{x_0  f_{13}f_{24}f_{12}f_{15}f_{25} }{p_{034}^3 p_{024} f_0^2 f_3^2 f_4}.
\end{align*}
Now $\det J(\xi_{05})\neq 0$  at a general point  in $Q_{05}\cap W$, so $Q_{05}$ is not contracted by $\psi_{05}$.

\underline{Case III: $\alpha=12$.}
By Corollary \ref{salt}, there exists some $M_{12}\in \PGL(4)$ such that 
\begin{align*}
	M_{12}\circ \psi=& [s''_0=ds_1'+cs_2:s_1:s_2:s''_3=r s_3'+t s_1]\\
=&[ p_{245} f_3 f_{13} f_{12}:  p_{045} f_3 f_4 f_{24}:   p_{234} f_0 f_3 f_{13}: p_{135}  f_4 f_{24}f_{12}].
\end{align*}
Therefore we prove the same statement for $M_{12}\circ \psi$, where $Q_{12}$ is mapped to the line $l_{12}=\{[0:*:*:0]\}$.
Then the same argument as in Case II reduces the proof to the Jacobian determinants of $\xi_{12}:\mathbb{A}^3\dashrightarrow \mathbb{A}^3$, where
\[\xi_{12}(x)= (y_0, y_2,b)=\left(\frac{s_0''}{s_1},\frac{s_2}{s_1},\frac{s_3''}{s_0''}\right).\]
Then Lemma \ref{imagesJac} and Proposition \ref{exc} show that the $\det J(\xi_{12})\neq 0$  at a general point  in $Q_{12}\cap W$. Hence $Q_{12}$ is not contracted by $\psi_{12}$.
\qed

\section{Images of the Exceptional Divisors}\label{imageE}
Here in this section we first show that the birational map $\phi:Y \dashrightarrow Y$ does not contract the following exceptional divisors: $E_0,E_3,E_4$ and $E_{12}, E_{15}, E_{25}$. To summarize, we show that $\phi$ maps $\{E_0,E_3,E_4\}$ birationally to $\{E_5,E_1,E_2\}$, and $\{E_{12}, E_{15}, E_{25}\}$ birationally to $\{E_{34}, E_{03}, E_{04}\}$. Then we prove that $\phi$ is a pseudo-automorphism of $Y$, and $\phi_X$ restricts to Keum's pseudo-automorphism $\kappa$.

\subsection{Image of $E_4$}
By symmetry of $E_0,E_3$ and $E_4$, we need only show the following:
\begin{prop}
	$\phi(E_4)\subseteq E_2$, and the restriction $\phi_{\mid E_4}:E_4\ra E_2$ is birational.
	\label{E4image}
\end{prop}
We first recall a lemma on linear systems of quartics in $\PP^2$ with base points. Let $\pi :W\ra \PP^2$ be the blow-up at six distinct points $a_1,\cdots, a_6$ such that (1) $a_1,a_2,a_3$ are not collinear; (2) $a_i$ is not on the three lines $\overline{a_1 a_2}, \overline{a_1 a_3}$ and $\overline{a_2 a_3}$ for $i=4,5,6$ and (3) no conic passes through all the six points. Consider the divisor class $C\sim 4h-2({e}_1+{e}_2+{e}_3)-({e}_4+{e}_5+{e}_6)$, where $h$ is the hyperplane class and ${e}_i$ is the exceptional divisor over $a_i$. Then we have
\begin{lemma}
	The complete linear system $\vv{C}$ has dimension $2$ and induces a birational morphism: $\alpha:W\ra \PP^2$.
	\label{del3}
\end{lemma}
\pf. First we show $\dim \vv{C}=2$. Identify $\PP^2\cong \proj\CC[x:y:z]$. We can assume $a_1=[1:0:0]$, $a_2=[0:1:0]$ and $a_3=[0:0:1]$. Since $a_4,a_5,a_6$ are not collinear with any two of $a_1,a_2,a_3$, we can assume $a_4=[1:1:1]$, $a_5=[1:u:v]$, and $a_6=[1:t:w]$, with $u,v,t$ and $w$ nonzero. Now any quartic polynomial $f$ vanishing at $a_1,a_2,a_3$ with multiplicity $2$ has the form $f=r_1 x^2 y^2+r_2 x^2 z^2+ r_3 y^2 z^2+ r_4 x^2 yz +r_5 xy^2 z+r_6 xyz^2$. Then vanishing at each of $a_4,a_5, a_6$ gives a linear condition on $r_i$, which together gives a $(3\times 6)$-matrix $M$ with columns indexed by $r_i$.

We observe that the $(3\times 3)$-minor of $M$ at the columns $(r_2,r_4,r_6)$ is nonzero. Otherwise there exists a nonzero vector $v:=(0,r_2,0,r_4,0,r_6)^T$ such that $Mv=0$. This implies that $f=r_2 x^2 z^2+ r_4 x^2 yz+r_6 xyz^2=xz (r_2 xz +r_4 xy +r_6 yz)$ vanishes at $a_1,a_2,a_3$ twice and $a_4,a_5,a_6$ once. Since $v$ and $w$ are nonzero, $xz$ does not vanish on $a_4,a_5$ and $a_6$.  Hence $r_2 xz +r_4 xy +r_6 yz$ is a conic through the six points, a contradiction. As a result, $M$ has rank $3$, which implies that $H^0(W,\oO_W(C))=6-3=3$. Hence $\dim \vv{C}=2$.

Let $\sigma_{pqr}$ be the standard Cremona transformation of $\PP^2$ centered at three noncollinear points $p,q,r$. Then $\sigma_{a_1 a_2 a_3}$ maps $a_4,a_5, a_6$ to three distinct points $b_4$, $b_5$ and $b_6$. Since no conic passes through all the six points, $b_4, b_5$ and $b_6$ are not collinear. Then we define $\beta:=\sigma_{b_4 b_5 b_6}\circ \sigma_{a_1 a_2 a_3}$, which is birational.
Direct calculation shows that $\beta^*\oO_{\PP^2}(1)=\oO_W(C)$. Therefore $\alpha=\beta$. Since $\sigma_{pqr}$ is resolved by blowing-up its center $p,q$ and $r$, we find $\alpha$ is a morphism. \qed

\pfof{Proposition \ref{E4image}}. Consider $\phi_{D}:Y\dashrightarrow \PP^3$, the birational map induced by $\vv{D}$.
Recall Definition \ref{sothers} and Corollary \ref{salt} that there are constants $c,d$ such that  \[s''_0:=f_3 f_{13}p_{245}f_{12}=cs_2+ds_1'=cs_2+d(s_1-s_0).\] Hence 
there exists some $M\in\PGL(4)$ sending $[s_0:s_1:s_2:s_3]$ to $[s_0:s_1:s''_0: s_3'] $ in $\PP^3$. Then $\vv{D}$ is generated by the following sections: 
\begin{align}
(\bar{s}_0 x_{E_{12}}x_{E_{15}}x_{E_{25}},\bar{s}_1 x_{E_{4}}x_{E_{12}},\bar{s}''_0 x_{E_{3}}x_{E_{4}},\bar{s}_3' x_{E_{4}}x_{E_{25}})
	\label{twistE4}
\end{align}
Now the last three sections all vanish on $E_4$. Hence under this choice of coordinates, $\phi_D$ sends $E_4$ to $[1:0:0:0]=q_2$. Next we blow up $q_2$ in the target. Define $\phi':Y\dashrightarrow \bl_{q_2} \PP^3$ as the lift of $\phi_D$. We abuse notation and write $E_2$ for the exceptional divisor over $q_2$ in $\bl_{q_2} \PP^3$. By restricting to some affine charts covering $E_4$ and $E_2$, we find $\phi'_{\mid{E_4}}$ sends a point $x\in E_4$ to the point $([1:0:0:0],[\bar{s}_1x_{E_{12}}:\bar{s}''_0 x_{E_{3}}:\bar{s}_3' x_{E_{25}}]_{\mid {E_4}})$. Recall Lemma \ref{E4q2} that $h^0(Y,\oO_Y(D-E_4))=3$. Hence we find $(\bar{s}_1x_{E_{12}},\bar{s}''_0 x_{E_{3}},\bar{s}_3' x_{E_{25}})$ span $H^0(Y,\oO_Y(D-E_4))$, and $\phi'_{\mid{E_4}}: E_4\dashrightarrow E_2$ is induced by the restriction of $\vv{D-E_4}$ to $E_4$.

 So we restrict $\vv{D-E_4}$ to $E_4$ and show it induces the same rational map $E_4\dashrightarrow E_2\cong\PP^2$ with the complete linear system in Lemma \ref{del3}. 
 We denote by $a_i$ the point where the proper transform of the line $l_{i4}$ meets $E_4$, and  $\ell_{ij}$ the intersection of the proper transform $\tilde{\Gamma}_{ij4}$ of the plane $\Gamma_{ij4}$ with $E_4$, for $i,j\neq 4$.  The exceptional divisor $E_4$ in $Y$ is isomorphic to the blow-up of $\PP^2$ at $a_0,a_2$ and $a_3$, so $\Pic(E_4)=\Z\{h, e_0,e_2,e_3\}$, with $e_i$ the exceptional divisor over $a_i$. The restriction map $r_{4}:\Pic(Y)\ra\Pic(E_4)$ is given by $H\mt 0$, $E_4\mt -h$, $E_{i4}\mt {e}_i$ for $i=0,2,3$, and all the else $E_i, E_{ij}$ to $0$. Now  Table \ref{e4r} shows the restrictions of those $Q_\alpha$ and $\tilde{\Gamma}_{ijk}$ appearing in (\ref{twistE4}).
\begin{table}[ht]
    \centering
\[\begin{array}{r|l|l}
		\xi\in \Pic(Y)	 & r_4(\xi) & \txt{The zeroes of the restricted section} \\
		Q_3 &  h-{e}_2-{e}_3 & \ell_{23} \\
		Q_4 &  2h-{e}_0-2{e}_2-{e}_3 & \ell_{02}+ \ell_{23} \\
		Q_0 &  h-{e}_0-{e}_2  &\ell_{02} \\
		Q_{24},Q_{12}, Q_{25} & 2h-{e}_0-{e}_2-{e}_3  &\txt{The proper transform of a conic through $a_0,a_2$ and $a_3$} \\
		Q_{13}  & 2h-2{e}_0-{e}_2-{e}_3 &\ell_{02}+ \ell_{03} \\
		\tilde{\Gamma}_{034} & h-{e}_0-e_3 & \ell_{03}\\
		\tilde{\Gamma}_{045} & h-{e}_0  & \ell_{05}\\
		\tilde{\Gamma}_{245} & h-{e}_2 &\ell_{25} \\
		\tilde{\Gamma}_{134} & h-{e}_3 &\ell_{13}
\end{array}\]\caption{Restriction to $E_4$}
    \label{e4r}
\end{table}

Let $c_{\alpha}$ be the conic $Q_\alpha\cap E_4$, for $\alpha=24,12$ or $15$. 
In the following, we abuse notations and write $\ell_{ij}$  (and $c_\alpha$) for the polynomials defining the sections $\ell_{ij}$ (and $c_\alpha$) in $\PP^2$ (identifying $E_4$ with the blow-up of $\PP^2$ at $a_0,a_2$ and $a_3$).
Then up to nonzero scalars, $\phi_{\mid E_4}':E_4\dashrightarrow E_2$ is the rational map defined by
\begin{align*}
[\bar{s}_1x_{E_{12}}:s''_0 x_{E_{3}}:\bar{s}_3' x_{E_{25}}]\mid_{E_4}&=
[\ell_{05}\ell_{23}^2\ell_{02}c_{24}:\ell_{25}\ell_{02}\ell_{03}\ell_{23}c_{12}:\ell_{13}\ell_{02}^2\ell_{23}c_{24}]\\
&=[\ell_{05}\ell_{23}c_{24}:\ell_{25}\ell_{03}c_{12}:\ell_{13}\ell_{02}c_{24}].
\end{align*}
Now we claim that $\phi'_{\mid E_4}$ is induced by the complete linear system of 
\[L\sim  4h-2({e}_0+{e}_2+{e}_3)-({e}_p+{e}_q+{e}_r)\] for some additional points $p,q,r$ in $\PP^2$ such that (1) $a_0,a_2,a_3$ are at linearly general position, (2) $p,q,r$ are not on the lines between $a_0,a_2$ and $a_3$, and (3) there is no conic through all the six points. Then by Lemma \ref{del3}, $\phi'_{\mid E_4}$ is birational, so $\phi_{\mid E_4}$ is birational, which finishes the proof.

We find the points $p,q$ and $r$ first. We make the following definition: 
\begin{itemize}
	\item Let $r$ be the unique point in $\ell_{05}\cap \ell_{13}$;
	\item let $p$ be the unique point in $\ell_{12}\cap c_{12}= \ell_{12}\cap c_{24}$ beside $a_0,a_2$ and $a_3$;
	\item let $q$ be the unique point in $\ell_{25}\cap c_{25}=\ell_{25}\cap c_{24}$ beside $a_0,a_2$ and $a_3$.
\end{itemize}
Indeed $r$ is well-defined. Here $p$ and $q$ are symmetric under the $\mathcal{S}_3$-action, so we need only show $p$ is well-defined. We show the explicit polynomials defining $c_{12}$ and $c_{24}$ and find $p$ as follows. Recall (\ref{f12}), (\ref{f24}), where we place $p_0,\cdots, p_5$ at $([0: 1: 0: 0], [1: 0: 0: 0], [1: a: b: c], [0: 0: 1: 0], [0: 0: 0: 1], [1: 1: 1: 1])$. Then

(i) $Q_{12}$ is defined by $f_{12}$  in (\ref{f12}) in $\PP^3$. Locally at $p_4=[0:0:0:1]$, we can take $U=\{ ((x,y,z),[X:Y:Z])\mid xY=yX, xZ=zX, yZ=zY\}$, and then identify $E_4$ with $\proj \CC[X:Y:Z]$.
Consider $f_{12}(x,y,z,w)$. Then on the affine chart $X=1$, we have $f_{12}(x,xY,xZ,1)=x^d h_{12}' (x,Y,Z)$ for some $d>0$ and polynomial $h_{12}'$ such that $x\nmid h_{12}' (x,Y,Z)$.
Thus on the affine chart $X=1$, $Q_{12}\cap E_4$ is defined by $h_{12}'(0,Y,Z)$. Homogenizing $h_{12}'(0,Y,Z)$ gives us a homogeneous polynomial $h_{12}(X,Y,Z)$ which define  $Q_{12}\cap E_4$. In our case we find
\[h_{12}=-a ( b-1)^2 c YZ 	-a b (1 - 2 c +   b c)   Y (X - Z)- a^2 b ( c-1)  X (X - Z).\]

(ii) Similarly, use $f_{24}$ in (\ref{f24}). Then  $Q_{24}\cap E_4$ is defined by \begin{align*}
h_{24}=&-(a - b) (a - c) (b - c)Y (X - Z)  - a (a - b) b (a - c) ( c-1) X (X - Z)+\\
&a (a - b) ( b-1) (c-1) c X Z+(a-1) (b-1) b (b - c) c X Y  -a (b-1)^2 (a - c) c X Y.
\end{align*}
Now $a_0=[0:1:0], a_2=[1:a:b]$ and $a_3=[0:0:1]$, and $\ell_{12}=\{bY-aZ=0\}$, so
\[c_{12}\cap \ell_{12}=c_{24}\cap \ell_{12}=\{a_0,a_2,a_3,[(b-c):a(1-c):b(1-c)]\}.\] Therefore 
\begin{align}
	\label{ypoint}
	p=[(b-c):a(1-c):b(1-c)].
\end{align}
is well defined. We also obtain that $r=[1:0:1]$. We next claim:
\begin{enumerate}
	\item $r\in c_{12}\cap c_{25}$;
	\item $p,q,r$ are distinct points, not on the lines $\ell_{02},\ell_{03}$ and $\ell_{23}$;
	\item $c_{24}, c_{12}$ and $c_{25}$ are smooth conics;
	\item No conic passes through all the six points $a_0,a_2,a_3,p,q,r$.
\end{enumerate}
Indeed, suppose $p=q$. Then $p\in  \ell_{12}\cap \ell_{25}$, so $p=a_2$, contradiction. The other claims in (1) and (2) follow from a direct calculation noticing the symmetry between $c_{12}$ and $c_{15}$, and between $p$ and $q$.
For (3), we need only to show each of the conic is irreducible. Equivalently, we can show the conics do not contain the three lines $l_{02}, l_{03}$ and $l_{23}$ by a calculation.
Finally for (4), suppose there is a conic $C$ through $a_0,a_2,a_3,p,q$ and $r$. Since $c_{24}$ is smooth, $c_{24}$ is uniquely determined by the five distinct points $a_0,a_2,a_3,p$ and $q$ on it. Therefore $C=c_{24}$ and $r\in c_{24}$, which contradicts the direct calculation that $r=[1:0:1]\not\in c_{24}$.

As a conclusion, each of the conics $c_{12}, c_{24}$ and $c_{25}$ passes through $a_0,a_2,a_3$ and exactly two of $p, q,r$. Therefore all the sections $\{\ell_{05}\ell_{23}c_{24},\ell_{25}\ell_{03}c_{12},\ell_{13}\ell_{02}c_{24}\}$ vanish at  $a_0,a_2,a_3,p,q$ and $r$ with multiplicities $(2,2,2,1,1,1)$. This proves that $\phi'_{\mid E_4}$ is the birational morphism induced by the complete linear system of $L$, which finishes the proof. \qed

\subsection{Image of $E_{12}$}
By symmetry of $E_{12},E_{15}$ and $E_{25}$, we need only show the following:
\begin{prop}
	\label{E12image}
	$\phi(E_{12})\subseteq E_{34}$. The restriction $\phi_{\mid E_{12}} : E_{12}\dashrightarrow E_{34}$ is induced by the complete linear system $\vv{\oO(1,1)}$ on $E_{12}\cong \PP^1 \times \PP^1$ and is birational.
\end{prop}
\pf. The result is local and we prove it without assuming that $\psi$ is birational. We restrict $\phi_D$ to $E_{12}\cong \PP^1\times \PP^1$. In $E_{12}$, we denote by $e$ the class of a section from $\overline{p_1 p_2}$, and $f$ the class of a fiber. Then $\Pic(E_{12})\cong \Z e+\Z f$. Under this identification, the restriction map $r_{12}$ sends $H,E_1$ and $E_2$ to $f$, $E_{12}$ to $-(e+f)$, and every other $E_i$ or $E_{ij}$ to the class $0$. As a result, $7$ out of the $9$ quartics $Q_\alpha$ restrict to $0$ except $Q_3$ and $Q_4$, where $r_{12}(Q_3)=r_{12}(Q_4)=e$. Finally, let $m_{ijk}=\vv{\{1,2\}\cap \{i,j,k\}}$. Then $r_{12}(\tilde{\Gamma}_{ijk})=e$ if $m_{ijk}=2$, $0$ if $m_{ijk}=1$ and $f$ if $m_{ijk}=0$. By Proposition \ref{wd} the map $\phi_{D}:Y\dashrightarrow \PP^3$ is given by 
$(\bar{s}_0 x_{E_{12}}x_{E_{15}}x_{E_{25}},\bar{s}_1 x_{E_{4}}x_{E_{12}},\bar{s}_2 x_{E_{3}}x_{E_{15}},\bar{s}_3 x_{E_{0}}x_{E_{25}})$. 
Therefore ${\phi_D}_{\mid E_{12}}$ is given by
\[[0:0:{\bar{s}_2}  x_{E_{3}}x_{E_{15}}\mid_{E_{12}}: {\bar{s}_3}  x_{E_{0}}x_{E_{15}}\mid_{E_{12}}].\]
This proves that $\phi_D(E_{12})\subset l_{34}$ by Definition-Theorem \ref{6dual}. Now up to scalars,
\[[{\bar{s}_2}  x_{E_{3}}x_{E_{15}}\mid_{E_{12}}: {\bar{s}_3}  x_{E_{0}}x_{E_{15}}\mid_{E_{12}}]=[e_3:e_4],\]
where $e_3:={Q_3}_{\mid E_{12}}$ and $e_4:={Q_4}_{\mid E_{12}}$. We see above that both $e_3, e_4\sim e$. We claim that $e_3\neq e_4$,  so that $\phi_D$ does not contract $E_{12}$ to a point. To see this, we restrict the sections to $S$. We find $E_{12}\cap S=T_{12}$. Then
\begin{align*}
e_3\cap S=Q_{3}\cap T_{12}&=\kappa\ik E_3\cap T_{12}=\kappa\ik (E_3\cap T_{34}).\\
	e_4\cap S=Q_{4}\cap T_{12}&=\kappa\ik E_{4}\cap T_{12}=\kappa\ik (E_4\cap T_{34}).
\end{align*}
Now $E_3\cap T_{34}\neq E_4\cap T_{34}$. Hence $e_3\neq e_4$. As a result, $\phi_D(E_{12})= l_{34}$.

Finally we blow up $l_{34}$ in the target to obtain the birational map $\phi_{12}:X\dashrightarrow \bl_{l_{34}}\PP^3$. Up to nonzero scalars, the restriction of $\phi_{12}$ to $E_{12}$ is given by
\begin{align*}
([0:0:e_3:e_4],[s_0 x_{E_{15}}x_{E_{25}},s_1 x_{E_{4}}]\mid_{E_{12}})=&([0:0:e_3:e_4],[\tilde{\Gamma}_{034}:\tilde{\Gamma}_{045}]\mid_{E_{12}})\\
=&([0:0:e_3:e_4],[\xi_{034}:\xi_{045}]).
\end{align*}
where $\xi_{034}:=(\tilde{\Gamma}_{034})_{\mid E_{12}}$ and $\xi_{045}:=(\tilde{\Gamma}_{045})_{\mid E_{12}}$. Then both $\xi_{034}, \xi_{045}\sim f$. Now in $\PP^3$, $\overline{p_1 p_2}\cap \Gamma_{034}\neq \overline{p_1 p_2}\cap \Gamma_{045}$, hence $\xi_{034}\neq \xi_{045}$. As a result,  ${\phi_{12}}_{\mid E_{12}}:E_{12}\dashrightarrow E_{34}$ is induced by   $\vv{\oO(1,1)}$, hence birational. Therefore the same results hold for $\phi_{\mid E_{12}}$.\qed

\subsection{$\phi$ is pseudo-automorphism}
\begin{theorem}
	For $p_0,\cdots, p_5$ in very general position, $\phi:Y\dashrightarrow Y$ is a pseudo-automorphism.
	\label{psdautoY}
\end{theorem}
\pf. By Theorem \ref{m}, $\psi$ is birational, and we can choose the coordinates of $p_i$ so that $\psi\ik=\psi'$. Therefore $\phi:Y\dashrightarrow Y$ is birational, whose inverse is the unique birational map $\phi'$ lifting $\psi'$. Applying Lemma \ref{psdhyp}, we need only show $\phi$ and $\phi\ik=\phi'$ do not contract any divisors. By Proposition \ref{exc}, $\psi$ only contracts the nine $Q_\alpha$. By Section \ref{imageQ}, $\phi$ does not contract the $9$ quartics $Q_\alpha$. Furthermore, $\phi$ is \'{e}tale at a general point $x$ in $Q_\alpha$ with $\phi(x)\in E_\alpha$. This shows that $\phi'$ is \'{e}tale at a general point of $E_\alpha$ in the target copy of $Y$. Now apply the symmetry between the linear system $D'$ and $D$ defining $\psi'$ and $\psi$ respectively. We find $\phi$ is \'{e}tale at a general point of $P_\beta$ in the source copy of $Y$. Therefore $\phi$ does not contract the $9$ exceptional divisors $E_{\beta}$ for $\beta\in \mathcal{B}$. Finally, the only divisors left are $E_0,E_3,E_4, E_{12},E_{15}$ and $E_{25}$, which $\phi$ does not contract by Propositions \ref{E4image} and \ref{E12image}. As a conclusion, $\phi$ contracts no divisors of $Y$. By symmetry, $\phi\ik$ contracts no divisors too. Hence $\phi$ is a pseudo-automorphism.\qed

Since $\phi$ is a pseudo-automorphism, $\phi$ will map any effective divisor of $Y$ birationally onto its image. In particular we conclude that $\phi:Q_\alpha\dashrightarrow E_\alpha$ is birational. This also proves that each $Q_\alpha$ is a rational quartic. See Remark \ref{qrmk}.

\subsection{Restriction of $\phi_X$ is Keum's automorphism}
We consider the images of the remaining $6$ lines under $\psi$. We abuse notation and write $l_{ij}$ for either $\overline{p_i p_j}$ or $\overline{q_i q_j}$. Define $\psi(l_{ij})$ as the closure of $\psi(U)$, such that  $U\subset l_{ij}$ is the open set where $\psi$ is defined. Then we have
\begin{prop}\label{6lines}
	\begin{align*}
	\psi(l_{14})=l_{02}, \quad \psi(l_{02})=l_{35}, \quad \psi(l_{35})=l_{14}, \quad	\psi(l_{01})=l_{45}, \quad \psi(l_{45})=l_{23}, \quad \psi(l_{23})=l_{01}.
	\end{align*}
\end{prop}
\pf. By symmetry we need only show $\psi(l_{14})=l_{02}$. We claim that this follows from that $\psi(\Gamma_{124})=\Gamma_{025}$ and $\psi(\Gamma_{134})=\Gamma_{024}$. Indeed, by Proposition \ref{exc}, the exceptional set of $\psi$ equals the $9$ quartics $Q_\alpha$. Since $l_{14}\not\subset Q_\alpha$ (Theorem \ref{pencilsQ}(3)),  $\psi$ is defined on an open subset of $l_{14}$ and does not contract $l_{14}$. Thus if our claim holds, then $\psi(l_{14})\subset \Gamma_{025}\cap \Gamma_{024}=l_{02}$. Then $\psi(l_{14})=l_{02}$ by the  irreducibility of $l_{04}$ and $l_{02}$.

So we prove the claims above. Since $\Gamma_{ijk}$ are irreducible and $\psi$ does not contract any $\Gamma_{ijk}$, we need only prove the inclusions. 
By Proposition \ref{planefusion}, $q_{025}(\psi)=p_{124}f_0 f_{05} f_{25}$. Therefore $\psi(\Gamma_{124})\subset \Gamma_{025}$. Next, by an argument similar to the proof of Proposition \ref{planefusion}, we find that $q_{024}(\psi)=p_{134}f_0 f_{4} f_{24}$, using that $f_0 f_4 f_{24}$ vanishes at $l_{13}$ and $l_{34}$ with multiplicities exactly $3$ and $2$, while $q_{024}(\psi)$ vanishes at $l_{13}$ and $l_{34}$ with multiplicities $4$ and $3$. Therefore  $\psi(\Gamma_{134})\subset \Gamma_{024}$.  \qed

\begin{corollary}\label{psdautoX}
The birational automorphism $\phi_X:X\dashrightarrow X$ is a pseudo-automorphism.
\end{corollary}
\pf. Since $\psi$ maps $l_{14}$ birationally to $l_{02}$,  the map $\phi_X:X\dashrightarrow X$ maps $E_{14}$ birationally to $E_{02}$. By symmetry, none of the exceptional divisors over the lines in Proposition \ref{6lines} are contracted by $\phi_X$. The same result for $\phi_X'$ by symmetry. Hence  $\phi_X$ is a pseudo-automorphism.\qed

Now we determine the pullback map of $\phi_X$ on $\Pic(X)$.
\begin{prop}\label{pullback}
	The pullback map $\eta:\Pic(X)\ra \Pic(X)$ induced by $\phi_X$ is given by:
	\[\begin{array}{llllll}
		\eta(H)=D,\\
		 	\eta(E_0)=Q_0, & 	\eta(E_3)=Q_3, &	\eta(E_4)=Q_4, 	&\eta(E_1)=E_3, & 	\eta(E_2)=E_4, &	\eta(E_5)=E_0, \\
			\eta(E_{05}) =Q_{05}, & \eta(E_{13}) =Q_{13}, &	\eta(E_{24}) =Q_{24}, &		\eta(E_{12}) =Q_{12}, & \eta(E_{15})=Q_{15}, &	\eta(E_{25}) =Q_{25}, \\	
\eta(E_{03}) =E_{15}, & \eta(E_{04}) =E_{25}, &	\eta(E_{34}) =E_{12},\\
		 \eta(E_{02})=E_{14}, & \eta(E_{35})=E_{02}, & \eta(E_{14})=E_{35}, &	 	 \eta(E_{45})=E_{01}, & \eta(E_{23})=E_{45}, & \eta(E_{01})=E_{23}.
	\end{array}\]
\end{prop}
\pf. We need only show $D$ has no fixed part, so that $\eta(H)=D$. Then the rest follows from Theorem \ref{psdautoY} and Proposition \ref{6lines}. By (\ref{Dsection}), the base locus of $D$ is supported on some of the pairwise intersections among $Q_\alpha$, $E_i$, $E_{ij}$ and $p_{ijk}$. Since those $Q_\alpha$ are distinct and irreducible by Theorem \ref{pencilsQ}, these intersections cannot contain any divisors. Hence $D$ has no fixed part. \qed

Now that we know the action of $\phi$ on $\Pic(X)$, we can use a computer program to verify that the matrix $M_\eta$ of $\eta$ has infinite order (for instance, compute the Jordan canonical form of $M_\eta$), so that $\phi$ has infinite order. Alternatively, we can inductively show the repeated images of $E_{i}$ and $E_{ij}$ under $\phi\ik$ span infinite many extremal rays in $\Effb(Y)$. We show the following examples.
\begin{theorem}
	Let $F_k:=(\phi\ik)^k(E_{03})$ for $k\geq 1$. Then
	\begin{align*}
F_k\sim\quad  & 2k(k-1)H-k(k-1)\sum_{i=0}^5 E_i\\
		 &-m_k(E_{03}+E_{04}+E_{34})-m_{k-1}(E_{12}+E_{15}+E_{25})-n_k(E_{05}+E_{13}+E_{24})-G_k,
	\end{align*}
	where	
\[\begin{array}{lll|l}
		G_k & m_k & n_k & k\\\hline
	E_{04}+E_{34}   & 3d^2-1     & 2(3d-1)d    &k=3d,\\
	E_{12}+E_{25}   & d(3d+2)    & 2(3d+1)d    &k=3d+1,\\
	-(E_{05}+E_{13})& (3d+1)(d+1)& 2(3d^2+3d+1)&k=3d+2.
\end{array}\]
In particular, every $F_k$ spans a different extremal ray in $\Effb(Y)$ and $\Effb(X)$. Thus $\phi$ and $\phi_X$ have infinite order in $\Bir(Y)$ and $\Bir(X)$. In particular $X$ and $Y$ are not Mori Dream Spaces.
	\label{inforder}
\end{theorem} 
\pf. The formula of $F_k$ follows from an induction on $k$. Indeed, let the right-hand side be $F'_k$. First, $ F'_1=E_{15}=\eta(E_{03})\sim F_1$.  Suppose $F_k=F'_k$. Then we can show
\begin{align*}
F_{k+1}-F'_{k+1}&=(6k^2-8k-6(n_k+m_{k-1}))(2H-\sum_{i=0}^5 E_i)\\
&-(4k(k-1)-4n_k-3m_{k-1}-m_{k+1})(E_{03}+E_{04}+E_{34}) \\
&-(4k(k-1)-3n_k-4m_{k-1}-n_{k+1})(E_{05}+E_{13}+E_{24})-\eta(G_k)+G_{k+1}.
\end{align*}
Then we can check that $6(n_k+m_{k-1})=6k^2-8k-4\alpha_k$, and $4n_k+3m_{k-1}+m_{k+1}=3n_k+4m_{k-1}+n_{k+1}=4k(k-1)-2\alpha_k$, with $\alpha_k=0$ if $k=3d$, $1$ if $k=3d+1$, and $-1$ if $k=3d+2$. Hence $F_{k+1}-F'_{k+1}=2\alpha_k A-\eta(G_k)+G_{k+1}=0$, where $A$ is defined in Definition \ref{DQ}. Therefore by induction $F_{k}=F'_{k}$ for all $k\geq 1$.

Now no pairs among those $F_k$ are scalar multiples of each other. Otherwise, say $F_a=\lambda F_b$. Then $a\e b \m 3$ since the two terms appearing in $G_k$ have different coefficients in $F_k$ from the third term.
Say $3\mid a$. Compare the coefficients of $E_{03}$ and $E_{04}$. Then \[\lambda=\frac{m_a}{m_b}=\frac{m_a+1}{m_b+1}.\]
Therefore $m_a=m_b$, which with $3\mid a, b$ implies that $a=b$. The cases when $a\e \pm 1 \m 3$ are the same.

Since $F_k$ span different extremal rays of $\Effb(Y)$, $\phi$ has infinite order. Now $\Effb(Y)$ and $\Effb(X)$ are not rational polyhedral, so $Y$ and $X$ are not Mori Dream Spaces.\qed

Now we return to $S$.
Let $O(\NS(S))^+$ be the group of isometries of $\NS(S)\cong \Pic(S)$ which leaves the set of effective divisors invariant. Recall that for a lattice $L$, the discriminant group of $L$ is the group ${L^*}/L$, which is finite abelian.  Let $D_S$ be the discriminant group of $\NS(S)$.
\begin{prop}\cite[Thm. 4.1]{Keum1997}
	\label{autS}
	Suppose $S$ is a Jacobian Kummer K3 surface with $\rho(S)=17$. Then
	\[\Aut(S)\cong \{f\in O(\NS(S))^+\mid f=\pm \id \txt{on } D_S\}.\]
\end{prop}

\begin{theorem}	The pseudo-automorphism $\phi_X:X\dashrightarrow X$ restricts to Keum's automorphism $\kappa:S\ra S$ associated to the Weber Hexad $\mathcal{H}=\{5,23,1,14,2,12\}$.
	\label{restrict}
\end{theorem}
\pf. By Proposition \ref{pullback} we know $\eta(S)=S$, so $\phi_X$ fixes the class of $S$. Since $S$ is the unique anticanonical section of $X$ (Lemma \ref{uniqueminusK}), we must have ${\phi_X}_{\mid S}:S\dashrightarrow S$ is a birational map. Now $S$ is K3, hence a minimal surface. Thus any birational automorphism of $S$ is in fact regular everywhere. Hence ${\phi_X}_{\mid S} \in \Aut(S)=\Bir(S)$. 

Now $D_S$ is uniquely determined by $\NS(S)$, therefore for any $g\in \Aut(S)$, the action of $g$ on $D_S$ is determined by the pullback $g^*:\NS(S)\ra \NS(S)$. Therefore by Proposition \ref{autS}, for $g, h\in \Aut(S)$, $g=h$ if and only if the induced linear maps of pullback $g^*=h^*$ on $\NS(S)$. Hence if we show the restriction of $\eta$ to $\NS(S)$ agrees with $\kappa^*$, then we must have ${\phi_X}_{\mid S}=\kappa\in \Aut(S)$, which proves the theorem.

Finally, the restriction map: $r:\Pic(X) \ra \Pic(S)$ satisfies $r(H)=H_S$, $r(E_i)=E_i$, and $r(E_{ij})=T_{ij}$. Proposition \ref{pullback} and Theorem \ref{pencilsQ} show that the restriction map $\eta$ to $\Pic(S)$ agree with $\kappa^*$ on the $\Q$-basis  $\{H_S, E_i, T_{ij}\}$ of $\Pic(S)$ (Proposition \ref{LH}). Hence $\eta$ restricts to $\kappa^*$, which finishes the proof.\qed

\begin{remark}
Indeed Theorem \ref{restrict} also implies that $\phi$ is of infinite order.
\end{remark}

\section{Cremona Transformations That Only Contracts Rational Hypersurfaces}\label{Cremona}
The linear system $\vv{D}$ in (\ref{Dclass}) induces the birational transformation $\psi:\PP^3\dashrightarrow \PP^3$ (Definition \ref{psi}).  Here we show $\psi$ contracts rational hypersurfaces only but is not generated by the standard Cremona transformation $\sigma_3$ and $\PGL(4)$. Recall that in \cite{Blanc2014}, the authors defined $G_n(\mathbf{k})$ as the subgroup of $\Bir(\PP^n)$ generated by $\sigma_n$ and $\PGL(n+1)$ over the field $\mathbf{k}$. They also defined $H_n$ to be the subset  of $\Bir(\PP^n)$ consisting of elements which contract rational hypersurfaces only. It is shown that $G_n\subseteq H_n$. On the other direction, the authors gave examples of birational transformations in odd dimensions that lie in $H_n$ but not $G_n$, hence showing $G_n\neq H_n$ when $n$ odd. In particular, they proved:
\begin{theorem}\cite[Thm. 1.4]{Blanc2014}
	Let $\mathbf{k}$ be any field and $n>2$ be odd. Suppose $H$ is an irreducible hypersurface which is sent by an element $g \in G_n(\mathbf{k})$ onto the exceptional divisor of an irreducible closed subset $Z$ (that is, the lift of $g$ to $\PP^3 \dashrightarrow \bl_Z \PP^3$ maps $H$ birationally onto $E_Z$). Then $Z$ has even dimension.
	\label{BH}
\end{theorem}
\begin{corollary}
	Let $\psi:\PP^3\dashrightarrow \PP^3$ be the birational map in Definition \ref{psi}. Then $\psi\in H_3$ but $\psi\not\in G_3(\CC)$.
	\label{notG3}
\end{corollary}
\pf. By Proposition \ref{exc}, $\psi$ only contracts the hypersurfaces $Q_\alpha'$. Each $Q_\alpha'$ is rational because they are birationally mapped to $E_\alpha$ which are rational (also see Remark \ref{qrmk}(4)). Hence $\psi\in H_3$. 
Now consider the rational hypersurface $Q_{12}'$. Here $Z:=l_{12}$ is a line, irreducible of dimension $1$, which is odd. By Theorem \ref{psdautoY} and Proposition \ref{QtoE}, $Q_{12}'$ is birationally mapped by $\psi$ onto $E_Z=E_{12}$, by Theorem \ref{BH}, $\psi\not\in G_3(\CC)$.\qed

\section{SQMs of Blow-ups of \texorpdfstring{$\PP^n$}{Pn} Along Points and Lines}\label{SQMpn}
We apply the construction by Castravet and Tevelev in \cite{Castravet2015} to construct an SQM for the blow-up of $\PP^n$ at $(n+3)$ general points and $9$ lines through $6$ of the $n+3$ points for $n\geq 3$. We define $Y_n$ to be the blow-up of $\PP^n$ at $(n+3)$ points at very general position and $9$ lines through six of them, such that when the six points we chose are indexed by  $\{0,1,2,3,4,5\}$, the $9$ lines are labeled by $(ij)\in\mathcal{I}=\{03,04,34,12,15,25,05,13,24\}$.
In particular, $Y_3=Y$. 
We prove the following:
\begin{prop}\label{YnSQM}
	For each $n\geq 4$ there is a small $\Q$-factorial modification (SQM) $\tilde{Y}_n$ of $Y_n$ such that $\tilde{Y}_n$ is a $\PP^1$-bundle over $Y_{n-1}$.
\end{prop}

\begin{corollary}\label{Yn}
	 For $n\geq 3$ and the $(n+3)$ points at very general position, $\Effb(Y_n)$ has infinitely many extremal rays, and $Y_n$ is not a Mori Dream Space.	
\end{corollary}
\pf. Suppose for some $n\geq 4$, $Y_n$ is a Mori Dream Space. Then the SQM $\tilde{Y}_n$ of $Y_n$ is a Mori Dream Space. By \cite{Okawa2016}, the surjection image $Y_{n-1}$ is also a Mori Dream Space. Inductively this proves that $Y_3=Y$ is a Mori Dream Space, which contradicts Theorem \ref{inforder}.

Now we prove that $\Effb(Y_n)$ has infinitely many extremal rays. Recall that a convex cone $\sigma$ is polyhedral if and only if $\sigma$ is spanned by finitely many rays, and equivalently, $\sigma$ has at most finitely many extremal rays. Therefore we only need to show $\Effb(Y_n)$ is not polyhedral. Suppose towards contradiction that $\Effb(Y_n)$ is polyhedral. Since the SQM $\tilde{Y}_n$ is isomorphic to $Y_n$ in codimension one, $\Effb(\tilde{Y}_n)\cong \Effb(Y_n)$. Hence $\Effb(\tilde{Y}_n)$ is polyhedral. By Lemma \ref{effP1bundle}, $\Effb(Y_{n-1})$ is polyhedral. Inductively we find $\Effb(Y_3)$ is polyhedral, which contradicts Theorem \ref{inforder}. \qed

\begin{lemma}
\label{effP1bundle}
Let $X$ be a normal projective variety. Let $p:\PP\ra X$ be a $\PP^1$-bundle. If $\Effb(\PP)$ is polyhedral, then so is $\Effb(X)$.
\end{lemma}
\pf. Here $\Pic(\PP)\cong p^*\Pic(X)\oplus \Z(\xi)$, where $\xi=\oO_{\PP}(1)$. Consider a divisor $D\in \Pic(X)$. We show $D\in \Effb(X)$ if and only if $p^*D\in \Effb(\PP)$ and $p^*D\cdot f=0$ for the fiber class $f$. 
Indeed, if $D$ is effective, then so is $p^*D$, and $p^* D\cdot f=0$. Conversely, since $p:\PP\ra X$ has connected fibers, we have $p_* \oO_\PP\cong \oO_X$. Hence 
\[H^0(\PP, p^*D)\cong H^0(X, p_*p^*D)\cong H^0(X,D\otimes p_* \oO_\PP)=H^0(X,D).\]
Hence if $p^*D$ is effective, then $D$ is effective. Then the claim follows from taking closures in $N^1(X)_\R$ and $N^1(\PP)_\R$. 
By the claim, $\Effb(X)$ equals the hyperplane section $f^\perp$ of $\Effb(\PP)$. Since $\Effb(\PP)$ is polyhedral, the hyperplane section $\Effb(X)$ is polyhedral. \qed

Recall Kapranov's blow-up construction of $\mzb{n}$ \cite{Kapranov1993} that $\mzb{n}$ is isomorphic to the successive blow-up of $\PP^{n-3}$ at $n-1$ points in linear general position, the lines, $2$-planes, $\cdots$, and all the linear subspaces of codimension at least $2$ through the $n-1$ points. Then the blow-up of $\mzb{n}$ at a very general point is a blow-up of $Y_{n-3}$ when $n\geq 7$. Now suppose the effective cone of the blow-up of $\mzb{n}$ at a very general point is polyhedral, then $\Effb(Y_{n-3})$ is also polyhedral, which contradicts Corollary \ref{Yn}. Thus we have proved:
\begin{corollary}
	For $n\geq 7$, the effective cone of the blow-up of $\mzb{n}$ at a very general point has infinitely many extremal rays. Hence the blow-up of $\mzb{n}$ at a very general point is not a Mori Dream Space.
	\label{mzbC}
\end{corollary}
\begin{remark}
	\label{10and6}
	We note that for $n\geq 10$, $\mzb{n}$ itself is not a Mori Dream Space, so the blow-up of $\mzb{n}$ at a very general point is not a Mori Dream Space. On the other hand, it is unknown whether the blow-up of $\mzb{6}$ at a general point is a Mori Dream Space.
\end{remark}

In the following we prove  Proposition \ref{YnSQM}. First, we review the definition of compatible sequences of sections (css) in \cite{Castravet2015}:
\begin{defi}\cite[Def. 5.2]{Castravet2015}
	Let $D_i$, $i=1,\cdots, N$ be irreducible divisors of a smooth variety $X$ with simple normal crossings. Call the pairwise and triple intersections among them $D_{ij}$ and $D_{ijk}$, and call the interior of them $D_{ij}^0$ and $D_{ijk}^0$. Assume that $D_{ij}^0$ and $D_{ijk}^0$ are either irreducible or empty. Fix $p:W\ra X$ a  $\PP^1$-bundle over $X$. Let $s_i:D_i\ra p\ik (D_i)$ be sections of $p$ over $D_i$, with images $Z_i$, for $i=1,\cdots, N$. Then we say the sections $s_i$ forms a compatible sequences of sections (css) if the following hold:
	\begin{enumerate}
		\item If $i<j$ and $D_{ij}\neq \emptyset$, then the following hold:
			\begin{enumerate}
				\item  $Z_i=Z_j$ over $D_{ij}$, or
				\item $Z_i$ and $Z_j$ are disjoint over $D^0_{ij}$ (interior of $D_{ij}$), in which case the locus in $D_{ij}$ where $Z_i$ and $Z_j$ agree is either empty or a union of subsets  $D_{ijk}$ for some indices $k$ such that $1\leq k< i$. Moreover for any such $k$ we have $Z_k=Z_i$ over $D_{ik}$; $Z_k=Z_j$ over $D_{jk}$; in addition, for any $p\in s_k(D^0_{ijk})$, the following relations between the tangent spaces hold
						\[T_{p,s_i(D_{ij})}\cap T_{p,s_j(D_{ij})} =T_{p,s_k(D_{ijk})}.\]
			\end{enumerate}
		\item If $D_{ijk}\neq \emptyset$, then there are $\{a,b\}\subset \{i,j,k\}$, $a\neq b$, such that $Z_a=Z_b$ over $D_{ab}$.
	\end{enumerate}\label{css}
\end{defi}

We refer to \cite{Maruyama1982} and \cite[Sec. 5]{Castravet2015} for discussions of elementary transformations of vector bundles. 
  Now recall
  \begin{prop}\cite[Prop. 5.4]{Castravet2015}\label{chain}
	  Given a css $s_i$ of the $\PP^1$-bundle $p:W\ra X$, with image $Z_i$, $i=1\cdots, N$, let $q:W^1\ra X$ be the elementary transformation of $p$ by the data $(D_1, Z_1)$. Then the proper transforms $Z'_i$ of $Z_i$ in $W^1$  for $i\geq 2$, form a css of $q$. Therefore iteratively there is a sequence of $\PP^1$-bundles $W^0=W$, $W^1$, $\cdots$, $W^N$ over $X$ such that $W^n$ is an elementary transformation of $W^{n-1}$.
\end{prop}

\pfof{Proposition \ref{YnSQM}}.
We first construct a css. We fix $(n+4)$ points $x, r_0,\cdots, r_{n+2}$ in $\PP^{n+1}$ in linearly general position, and consider the natural morphism $\pi: \bl_x \PP^{n+1}\ra \PP^{n}$ resolving the projection $\PP^{n+1}\dashrightarrow \PP^n$ from $x$. Then $\pi$ is a $\PP^1$-bundle over $\PP^n$. 
Let $p_i=\pi(r_i)$. We denote by $l_\alpha$ the linear subspace (a line or a point) passing through the points $\{p_i\mid i\in \alpha\}$, where $\alpha\in \{0,1,2,\cdots,n+2\}\cup \mathcal{I}$. 
Let $u_n:Y_n\ra \PP^n$ be the successive blow-up of $\PP^n$ at the $(n+3)$ points $r_0,\cdots,r_{n+2}$ and the $9$ lines indexed by $\mathcal{I}$, that is, all those linear subspaces $l_\alpha$. Let $\pi':W\ra Y_{n}$ be the pullback of the bundle $\pi$. 
As in \cite[Proof of 1.1]{Castravet2015}, we choose sections $t_\alpha: l_\alpha \ra \pi\ik(l_\alpha)$ such that the image of $l_\alpha$ is the linear subspace $L_\alpha$ passing through the corresponding points  $\{r_i\mid i\in \alpha\}$.  We call $D_\alpha$ the exceptional divisor in $Y_n$ over $l_\alpha$, and pull back $t_\alpha$ to a section $s_\alpha: D_\alpha\ra(\pi')\ik(D_\alpha)$. Call $Z_\alpha$ the image $s_\alpha(D_\alpha)$.

Next we check that those sections $s_i$ form a css under the increasing order 
\[\{0,1,2,\cdots,n+2, 03,04,34,12,15,25,05,13,24\}.\]
Indeed, if $i,j,k,l\in \{0,1,\cdots, n+2\}$ are distinct indices, then $D_i\cap D_j=\emptyset$, $D_{ij}\cap D_{kl}=D_{ij}\cap D_{ik}=\emptyset$ and $D_{i,ij}:=D_i\cap D_{ij}\cong\PP^1$ is a fiber, whenever the divisor $D_{ij}$ is defined. Thus every triple intersection among $D_\alpha$ is empty, making (2) of Definition \ref{css} true. Furthermore, this implies that $\{D_\alpha\}$ are indeed simple normal crossing. Finally, for $D_i$ and $D_{ij}$, we find $Z_i$ agrees with $Z_{ij}$ over $D_{i,ij}$ since they are pullbacks of the sections $t_i$ and $t_{ij}$ which agree over the point $p_i$. This proves that $\{D_\alpha\}$ form a css.

Applying Proposition \ref{chain}, there exists a chain of $\PP^1$-bundles $W_\alpha$ over $Y_n$ so that each is an elementary transformation of the previous one:
\[W, W_0,\cdots, W_{n+2},W_{03},W_{04},\cdots, W_{24},\] where $W_\alpha$ is the successive elementary transformation of $W$ about the data $(D_\beta,Z_\beta)$ for every $\beta\leq \alpha$.
On the other hand, we identify $Y_{n+1}$ with the blow-up of $\bl_x \PP^{n+1}$ at the $n+3$ points $r_i$ and the $9$ lines through $r_i$ indexed by $\mathcal{I}$. Denote by $v_\alpha:X_\alpha\ra \PP^{n+1}$  the intermediate blow-ups at $x$ and all the linear subspaces $L_\beta$ for $\beta\leq \alpha$. Then we have a chain of blow-ups:
\[\bl_x \PP^{n+1}, X_0,\cdots, X_{n+2},X_{03},X_{04},\cdots, X_{24}=Y_{n+1}.\]

Finally we show that $W_{24}$ is an SQM of $Y_{n+1}$, so that $W_{24}$ is the SQM $\tilde{Y}_{n+1}$ needed.
We are adopting the proof of Claim 3 in \cite[Proof of Thm. 1.1]{Castravet2015}, where the author proved the same result for $\mzb{n}$ as the successive blow-up of $\PP^{n-3}$. Here we blow up one extra point and only $9$ lines, which does not affect the original argument. 

For the reader's convenience, we recall their proofs here.
For each $W_\alpha$ above, there exists a birational map $\varphi_\alpha: W_\alpha \dashrightarrow X_\alpha$ which on an open locus coincides with the natural morphism $\varphi:W\ra \bl_x \PP^{n+1}$, which is identity on an open locus. 
First, since any elementary transformation of vector bundles keeps the Picard rank, we have $\rho(W_{24})=\rho(W)=1+\rho(Y_n)=\rho(Y_{n+1})$.  Then we only need to inductively show that each $\varphi_\alpha: W_\alpha \dashrightarrow X_\alpha$ is a birational contraction, in the sense that $\varphi_\alpha\ik$ contracts no divisors. If so, then $\varphi_{24}$ must be a small modification of $Y_{n+1}$ because $W_{24}$ and $Y_{n+1}$ have the same Picard rank, which proves the claim.

So we run an induction for $\alpha$.
The base case is clear: $\varphi:W\ra \bl_x \PP^{n+1}$ is a birational contraction since $u_n\ik$ contracts no divisors. Suppose $\varphi_\beta$ is a birational contraction and suppose $\alpha$ is the next index after $\beta$. 

Let $f:A_\alpha \ra W_\beta$ be the blow-up of $W_\beta$ along $Z_\alpha$, with exceptional divisor $G_\alpha$. Let $F_\alpha:=(\pi')\ik (D_\alpha)$. Then the elementary transformation $A_\alpha$ is given by $g:A_\alpha \ra W_\alpha$, the blow-down of the proper transform $\tilde{F}_\alpha$ of $F_\alpha$.
Using \cite[Lem. 5.6]{Castravet2015}, we only need to show
\begin{enumerate}
	\item The rational map $\tilde{\varphi_\alpha}:= \varphi_\alpha\circ g: A_\alpha\dashrightarrow X_\alpha$ is a birational contraction, and 
	\item The bundle $F_\alpha$ is contracted by $\varphi_\beta$ (so that $\tilde{F}_\alpha$ is contracted by $\tilde{\varphi}_\alpha$).
\end{enumerate}
The key observation is that we only need to prove both arguments on an open set $V_\alpha\subset W_\alpha$ which intersects $F_\alpha$. To this aim, we define $U_\alpha$ to be the complement in $\PP^n$ of all the linear subspaces $l_{\alpha'}\subset \PP^n$ of $\alpha'\lneq\alpha$. Then $V_\alpha:={\pi'}\ik u_n\ik (U_\alpha)$ (the scheme-theoretical preimage) is open and intersects $F_\alpha$. 
The elementary transformation construction commutes with base change, so when restricted to $V_\alpha$, we have a chain of elementary transformations of $\PP^1$-bundles over $ u_n\ik (U_\alpha)$, and the induced proper birational morphism
\[(\varphi_\alpha)_{V}:V_\alpha\ra v_\alpha\ik (U_\alpha).\]
Each elementary transformation $W_{\alpha'}$ with $\alpha'\lneq \alpha$, is an isomorphism outside its center $Z_{\alpha'}$. In particular, the restrictions of $W_{\alpha'}$ for $\alpha'\lneq \alpha$, and $W$ to $u_n\ik(U_\alpha)$ are all isomorphic. On the other hand, The blow-ups $X_{\alpha'}$ for $\alpha'\lneq \alpha$ are isomorphic over $U_\alpha$. As long as the six points $p_i$ are in linearly general position, $(\varphi_\alpha)_V$ is a local isomorphism at the generic point of the exceptional divisors $G_\alpha$ of $A_\alpha$, which is sent to $E_\alpha$ of $X_\alpha$. Hence by \cite[Lem. 5.5]{Castravet2015}, $\tilde{\varphi}_\alpha\ik$ does not contract $E_\alpha$. By the induction hypothesis, $\varphi_\beta$ is a birational contraction. Together we know no divisors in $X_\alpha$ are contracted by $\tilde{\varphi}_\alpha$. Therefore $\tilde{\varphi}_\alpha$ is a birational contraction, which proves (1).
When restricted to $V_\alpha$, the $\PP^1$-bundle $F_\alpha\cap V_\alpha\ra D_\alpha \cap u_n\ik(U_\alpha)$ coincides with the $\PP^1$-bundle $\pi'_{\mid D_\alpha \cap u_n\ik( U_\alpha)}:{\pi'}\ik (D_\alpha)\cap V_\alpha\ra D_\alpha \cap u_n\ik(U_\alpha)$, which is contracted by $\varphi_\beta$ to $\pi\ik(l_\alpha\cap U_\alpha)$. This proves (2).

As a conclusion, both claims above hold over $V_\alpha$, hence hold for $W_\alpha$. This proves that $\varphi_\alpha$ is a birational contraction. By induction the last one $\varphi_{24}$ is a birational contraction, hence proving $W_{24}$ is an SQM of $Y_{n+1}$.\qed

\bibliography{mybib.bib}
\bibliographystyle{halpha}

\end{document}